\documentclass [11pt]{amsart}
\usepackage {amssymb, amscd, url, color}
\usepackage[dvips,final]{graphicx}
\usepackage[text={6.5in,9in},centering,letterpaper,dvips]{geometry}
\usepackage{hyphenat}
\usepackage[dvips]{xy}

\newtheorem {theorem}{Theorem}[section]
\newtheorem {lemma}[theorem]{Lemma}
\newtheorem {prop}[theorem]{Proposition}
\newtheorem {corollary}[theorem]{Corollary}

\theoremstyle{definition}
\newtheorem {definition}[theorem]{Definition}

\theoremstyle{remark}
\newtheorem {remark}[theorem]{Remark}
\newtheorem {example}[theorem]{Example}

\newcommand\oLink{\vec{L}}
\newcommand\oKnot{\vec{K}}

\newcommand{\abs}[1]{\lvert#1\rvert}
\let\savecolon\colon
\renewcommand{\colon}{\mskip0.5mu\savecolon\thinspace}   
\renewcommand\th{^{\text{th}}}
\newcommand\Z{\mathbb{Z}}
\newcommand\R{\mathbb{R}}

\newcommand\Xs{\mathbb{X}}
\newcommand\Os{\mathbb{O}}
\newcommand\Pent{\mathrm{Pent}}

\newcommand\Rect{\mathrm{Rect}}
\newcommand\Ring{\mathcal R}

\def\HFa{\widehat {\mathit{HF}}}

\def\Interior{\mathrm{Int}}

\def\Cm{\mathit{C}^-}
\def\Ca{\widehat{\mathit{C}}}
\def\CKm{\mathit{CK}^-}
\def\CLm{\mathit{CL}^-}
\def\HFm{\mathit{HF}^-}

\def\CFa{\widehat{\mathit{CF}}}

\def\HFKm{{\mathit{HFK}}^-}

\def\HFLm{\mathit{HFL}^-}
\def\OneHalf{{\textstyle\frac{1}{2}}}

\def\cm{\cdot}
\def\s{\mathbf s}
\def\x{\mathbf x}
\def\y{\mathbf y}
\def\z{\mathbf z}
\def\Torus{T^2}
\newcommand\orL{\vec{L}}
\def\S{\mathbf S}

\def\NEoff{\mathcal I}
\def\NESW{\mathcal J}
\def\LegK{\mathcal K}
\def\oLegK{\vec{\mathcal K}}

\def\oLegL{\vec{\mathcal L}}
\DeclareMathOperator\lk{lk}
\def\deg{{}^\circ}

\def\TransK{\mathcal T}

\DeclareMathOperator{\writhe}{wr}

\def\Field{\mathbb F_2}
\newcommand\Filt{\mathcal F}

\newcommand\dm{\partial^-}

\newcommand\EmptyRect{\Rect^o}

\DeclareMathOperator{\TB}{tb}
\DeclareMathOperator{\Rotation}{r}
\DeclareMathOperator{\SL}{sl}

\newcommand\Left{\mathcal L}
\newcommand\Right{\mathcal R}
\newcommand\FL{F^L}
\newcommand\FR{F^R}
\def\pL{\pi^L}
\def\pR{\pi^R}

\def\Interval{\mathbf I}

\def\EmptyRect{\Rect^\circ}
\def\EmptyPent{\Pent^\circ}

\newcommand{\NW}{\text{\sl NW}}
\newcommand{\NE}{\text{\sl NE}}
\newcommand{\SW}{\text{\sl SW}}
\newcommand{\SE}{\text{\sl SE}}
\newcommand{\ONW}{\text{\sl O:NW}}
\newcommand{\ONE}{\text{\sl O:NE}}
\newcommand{\OSW}{\text{\sl O:SW}}
\newcommand{\OSE}{\text{\sl O:SE}}
\newcommand{\XNW}{\text{\sl X:NW}}
\newcommand{\XNE}{\text{\sl X:NE}}
\newcommand{\XSW}{\text{\sl X:SW}}
\newcommand{\XSE}{\text{\sl X:SE}}

\newcommand\zUR{\mathbf{z^+}}
\newcommand\zLL{\mathbf{z^-}}
\newcommand\xUR{\mathbf{x^+}}
\newcommand\xLL{\mathbf{x^-}}
\def\UR{\lambda^+}
\def\LL{\lambda^-}
\def\TransElt{\theta}

\def\Mirror{m}

\newcommand{\sign}{\mathcal S}

\def\graph#1{\includegraphics{draws/#1}}
\def\mathcenter#1{\vcenter{\hbox{$#1$}}}

\def\mfigb#1{\mathcenter{\includegraphics[trim=-2 -2 -2 -2]{draws/#1}}}

\begin{document}

\title{Legendrian knots, transverse knots and combinatorial Floer
homology}

\author[Ozsv\'ath]{Peter Ozsv\'ath}
\thanks {PSO was supported by NSF grant number DMS-0505811 and FRG-0244663}
\address {Department of Mathematics, Columbia University\\ New York, NY 10027}
\email {petero@math.columbia.edu}

\author[Szab{\'o}]{Zolt{\'a}n Szab{\'o}}
\thanks{ZSz was supported by NSF grant number DMS-0406155 and FRG-0244663}
\address{Department of Mathematics, Princeton University\\ Princeton, New Jersey 08544}
\email {szabo@math.princeton.edu}

\author[Thurston]{Dylan Thurston}
\thanks {DPT was supported by a Sloan Research Fellowship.}
\address {Department of Mathematics, Barnard College, Columbia University\\ New York, NY 10027}
\email {dthurston@barnard.edu}

\begin {abstract} 
  Using the combinatorial approach to knot Floer homology, we define
  an invariant for Legendrian knots (or links) in the three-sphere,
  with values in knot Floer homology. This invariant can also
  be used to construct an invariant of transverse knots.
\end {abstract}

\maketitle
\section{Introduction}

We use link Floer homology to study contact phenomena for links in
the three-sphere, endowed with its standard contact structure.
 
Knot (or link) Floer homology is an invariant defined using Heegaard
diagrams and holomorphic disks~\cite{Knots,RasmussenThesis}. It comes
in various forms, but the version which will be of primary interest to
us here is the variant which associates to a knot $K$ a module over
$\Z[U]$, denoted
\[
\HFKm(K)=\bigoplus_{m,a\in\Z}\HFKm_m(K,a),
\]
where $U$ acts as an endomorphism which is
homogeneous of degree~$-2$ for the Maslov grading~$m$ and degree~$-1$ for
the Alexander grading~$a$.

Manolescu, Ozsv\'ath, and Sarkar gave~\cite{MOS} an explicit
description of knot Floer homology for a knot in the three-sphere as
the homology groups of a chain complex $\CKm$ which is described in
terms of the combinatorics of a grid diagram for a knot.  In fact, the
constructions of~\cite{MOS} are done with coefficients in $\Z/2\Z$; a
lift of these constructions to coefficients in $\Z$ is given
in~\cite{MOST}, along with a purely combinatorial proof of the fact
that their homology groups are knot invariants, which entirely
circumvents the holomorphic description. Given a grid diagram~$G$ for
the mirror~$\Mirror(K)$ of a knot $K$, we refer to the resulting
complex as the
\emph{combinatorial chain complex} for~$G$, denoted $\CKm(G)$.
Specifically, $\CKm(G)$ is generated by a collection of generators
$\S$, equipped with grading functions $A$ and $M$, so that
\[
\CKm(G) = \bigoplus_{m,a \in\Z}\CKm_m(G,a)
\]
where $\CKm_m(G,a)$ denotes the set of
$\Z$-linear combinations of generators $\x$ with $A(\x)=a$ and
$M(\x)=m$.  $\CKm(G)$ has an explicitly
defined differential $\dm\colon \CKm(G)\longrightarrow \CKm(G)$,
with $\dm\colon \CKm_{m}(G,a)\longrightarrow \CKm_{m-1}(G,a)$
(so the bigrading descends to homology). On the other hand, a grid
diagram may also be used to give a combinatorial presentation of a
Legendrian knot~$\LegK$ or a transverse knot~$\TransK$; we will
exploit this fact here.  (By our conventions, explained in
Section~\ref{sec:GridMoves},
the topological type of~$\LegK$ is the mirror of~$K$).

Given a grid diagram for a knot, we will exhibit a concrete pair of
generators $\zUR$ and~$\zLL$ for the combinatorial chain complex, both
of which are cycles.  These elements are defined in
Definitions~\ref{def:Canonical} and~\ref{def:CanonicalWithSigns}, and
the fact they are cycles is established in Lemma~\ref{lemma:IsACycle}.
Our aim here is to study this pair of cycles. We show that the pair of
induced homology classes is an invariant for Legendrian knots. To
describe the bigradings of these elements, we use the two classical
invariants of a Legendrian knot~$\LegK$, the
\emph{Thurston-Bennequin invariant}~$\TB(\LegK)$ and the
\emph{rotation number}~$\Rotation(\oLegK)$,
which we recall in Section~\ref{sec:Contact}. The overall sign of the
rotation number depends on the orientation of~$\oLegK$.  (Note that
we restrict our attention in this introduction, and indeed through
most of the present paper, to the case of knots, as opposed to links,
though most of the results here carry over with minor modifications to
the case of links. These generalizations are discussed in
Section~\ref{sec:Links}.)

\begin{theorem} 
        \label{thm:LegendrianInvariant}
        For a grid diagram~$G$ which represents a knot, let
        $\oLegK=\oLegK(G)$ be the
        corresponding oriented Legendrian knot. Then there are two
        associated cycles $\zUR=\zUR(G)$ and $\zLL=\zLL(G)$, supported
        in bigradings
\begin{align*}
M(\zUR)&=\TB(\oLegK)-\Rotation(\oLegK)+1,&
M(\zLL)&= \TB(\oLegK)+ \Rotation(\oLegK) + 1,\\
A(\zUR)&=\frac{\TB(\oLegK)-\Rotation(\oLegK)+1}{2},& 
A(\zLL)&=\frac{\TB(\oLegK)+\Rotation(\oLegK)+1}{2}.
\end{align*}
Moreover, if $G$ and $G'$ are two different grid diagrams which
represent Legendrian isotopic oriented knots, then there is 
a quasi-isomorphism of chain complexes
$$\Phi\colon \CKm(G) \longrightarrow \CKm(G')$$
with
\begin{align*}
\Phi(\zUR(G))&=\zUR(G'),& \Phi(\zLL(G))&=\zLL(G').
\end{align*}
\end{theorem}

If $G$ is a grid diagram representing an oriented Legendrian knot
$\oLegK$, then we denote the homology classes of $\zUR(G)$ and $\zLL(G)$
in $\HFKm(\Mirror(K))$ by $\UR(\oLegK)$ and $\LL(\oLegK)$ respectively,
and refer to them as the \emph{Legendrian
invariants} of $\oLegK$. 
The pair of invariants $\UR(\oLegK)$ and $\LL(\oLegK)$ can be used to
distinguish Legendrian knots with the same classical invariants. For
example, for the two different Legendrianizations of the knot $5_2$,
$\oLegK_1$ and $\oLegK_2$ with $\Rotation=0$ and $\TB=1$, we have that
$\UR(\oLegK_1)\neq \LL(\oLegK_1)$, while $\UR(\oLegK_2)=\LL(\oLegK_2)$.
See Example~\ref{ex:Chekanov} below.

There is some symmetry in the construction of the cycles $\zUR$ and
$\zLL$. To this end, recall the \emph{Legendrian mirror} construction
\cite{Ng,FuchsTabachnikov}. Given an oriented Legendrian knot
$\oLegK$, one can rotate the front projection by $180^\circ$ around
the $x$-axis to obtain
the Legendrian knot projection of an oriented Legendrian knot denoted
$\mu(\oLegK)$. Classical invariants are related by
\begin{align*}
\TB(\mu(\oLegK))&= \TB(\oLegK), & 
\Rotation(\mu(\oLegK))&= -\Rotation(\oLegK).
\intertext{One can instead reverse the orientation of
  $\oLegK$ to obtain a different oriented Legendrian knot $-\oLegK$ with}
\TB(-\oLegK)&= \TB(\oLegK),&
\Rotation(-\oLegK)&= -\Rotation(\oLegK).
\end{align*}

\begin{prop}
  \label{prop:Symmetries}
  Suppose that $G$ is a grid representation of $\oLegK$, $G_1$ is a grid
  representation of the orientation reversal $-\oLegK$, and $G_2$ is a grid
  representation of the Legendrian mirror $\mu(\oLegK)$. Then there
  are quasi-isomorphisms
  \begin{align*}
    \Phi_1\colon \CKm(G) & \longrightarrow \CKm(G_1), &
    \Phi_2\colon \CKm(G) & \longrightarrow \CKm(G_2)
  \end{align*}
  which have the property that
  \begin{align*}
    \Phi_1(\zUR(G))&=\zLL(G_1),& \Phi_2(\zUR(G))&=\zLL(G_2),\\
    \Phi_1(\zLL(G))&=\zUR(G_1),& \Phi_2(\zLL(G))&=\zUR(G_2).
  \end{align*}
\end{prop}

More interestingly, this invariant behaves in a controlled manner
under stabilizations of the Legendrian knot. Specifically, recall that
one can locally introduce a pair of cusps in the front projection of a
Legendrian knot~$\LegK$ to obtain a new Legendrian knot~$\LegK'$
which is in the same topological type as~$\LegK$, but which is not
Legendrian isotopic to the original knot. The knot~$\LegK'$ is called
a \emph{stabilization} of $\LegK$. If we fix an orientation for
$\LegK$, we can distinguish the two ways of stabilizing as
\emph{positive} and
\emph{negative}.  (We adhere to the conventions spelled out, for
example, in Etnyre's survey~\cite{Etnyre}, which, incidentally, is
also an excellent
reference for the basic theory of Legendrian and transverse knots. We
review these conventions in Section~\ref{sec:Contact}.)
The Legendrian invariants transforms in the following manner under
stabilizations:

\begin{theorem}
        \label{thm:StabilizationTheorem} Let $\oLegK$ be an oriented
        Legendrian knot, and $\oLegK^-$ (respectively $\oLegK^+$) be the
        oriented Legendrian knots obtained as a single negative
        (respectively positive) stabilization of $\oLegK$. Then there are 
        quasi-isomorphisms
        \begin{align*} 
                \Phi^-\colon \CKm(\oLegK) &\longrightarrow \CKm(\oLegK^-), &
                \Phi^+\colon \CKm(\oLegK) &\longrightarrow \CKm(\oLegK^+) 
        \end{align*} 
        under which
        \begin{align*}
          \phi^-(\UR(\oLegK))&=\UR(\oLegK^-),&
              U\cm \phi^+(\UR(\oLegK))&=\UR(\oLegK^+),\\
          U\cm \phi^-(\LL(\oLegK))&= \LL(\oLegK^-),  &
              \phi^+(\LL(\oLegK))&= \LL(\oLegK^+),
        \end{align*}
        where $\phi^\pm$ denotes the map induced on homology by $\Phi^\pm$.
\end{theorem}

The above theorem suggests an application to transverse knots. A
transverse knot $\TransK$ inherits its orientation from the contact
structure on $S^3$. If we let $\oLegK$ be a Legendrian approximation
to $\TransK$ (cf.\ Section~\ref{sec:Contact}) with the natural induced
orientation, then we can define the \emph{transverse invariant}
$\TransElt(\TransK)$ to be $\UR(\oLegK)$.  By a result of Epstein,
Fuchs, and Meyer~\cite{TransverseIsotopic} (generalized by Etnyre and
Honda~\cite{EHTransverse}), any two Legendrian approximations to some
given transverse knot are equal after some number of negative
stabilizations. Thus, in view of
Theorem~\ref{thm:StabilizationTheorem}, we can conclude that the
transverse invariant is independent of the choice of Legendrian
approximation $\oLegK$ used in its definition.

\begin{corollary}
  \label{cor:TransverseWellDefined}
  The transverse invariant $\TransElt(\TransK)$ depends only on the
  transverse isotopy class of the transverse knot $\TransK$; i.e.,
  if $G$ and $G'$ are two grid diagrams representing two Legendrian approximations
  to $\TransK$, then there is a quasi-isomorphism
  $$\Phi\colon \CKm(G) \longrightarrow \CKm(G')$$
  whose induced map $\phi$ on homology has the property that
  $\phi(\TransElt(\TransK))=\TransElt(\TransK')$.
\end{corollary}

In later work by Ng, Ozsv\'ath, and Thurston, the invariant $\theta$ is used to
distinguish particular transversally non-isotopic knots with the same
classical invariants~\cite{NOT}.

These invariants also satisfy the following non-vanishing property.

\begin{theorem}
  \label{thm:NonVanishing} For any Legendrian knot $\oLegK$, the
  homology classes $\UR(\oLegK)$ and $\LL(\oLegK)$ are non-trivial; and
  they are not $U$-torsion classes (i.e., for all positive integers $d$,
  $U^d$ times these classes is non-trivial). Similarly, for any transverse
  knot~$\TransK$, the transverse invariant~$\TransElt(\TransK)$ 
  does not vanish, and indeed is not $U$-torsion.
\end{theorem}

We can use this to reprove bounds on the Thurston-Bennequin invariant.
Recall that $\tau(K)\in\Z$ is a knot concordance invariant defined
using $\HFKm(K)$~\cite{4BallGenus}.  One definition is that $\tau(\Mirror(K))$ is the
minimal
Alexander grading of any element $\xi$ of $\HFKm(K)$ for which
$U^d\xi\neq 0$ for all integers $d\geq 0$. (This is not the
usual definition of~$\tau$, but we verify in Appendix~\ref{sec:OnTau}
that it is equivalent.) Combining
Theorem~\ref{thm:NonVanishing} with
Theorem~\ref{thm:LegendrianInvariant}, we see at once that for any
Legendrian knot~$\LegK$,
\begin{equation}
        \label{eq:TBtau}
\abs{\Rotation(\LegK)}+\TB(\LegK)\leq
2\tau(\LegK)-1,
\end{equation}
a bound which was first proved~\cite{PlamenevskayaTB} using the
contact invariant in Heegaard Floer homology~\cite{HolDiskContact}.
The inequality
$\tau(K)\leq g(K)$ is easy to establish~\cite{Knots}, giving
another proof of Bennequin's inequality.  Indeed, since $\tau(K)\leq
g^*(K)$~\cite{4BallGenus},
we have yet another proof of the ``slice-Bennequin inequality''
\[
\abs{\Rotation(\LegK)}+\TB(\LegK)\leq 2g^*(K)-1,
\]
first proved using methods of gauge theory~\cite{Rudolph,KMmilnor}.

It is interesting to compare our results here with those of
Plamenevskaya~\cite{PlamenevskayaTransverse}. In that paper,
Plamenevskaya uses braid representatives to give an invariant for
transverse knots which takes values in the Khovanov homology of the
knot. She uses this to prove the following analogue of Equation~\eqref{eq:TBtau}:
$$\abs{\Rotation(\LegK)}+\TB(\LegK)\leq s(K)-1,$$
where now $s(K)$ is the Rasmussen invariant coming from
Khovanov homology~\cite{RasmussenSlice}. This bound, when combined
with Rasmussen's bound $s(K)\leq 2g^*(K)$, 
gives a different proof of the slice-Bennequin inequality. 
(Ng gives a different Legendrian Thurston-Bennequin bound using
Khovanov homology~\cite{NgLegTB}.)

In Section~\ref{sec:Contact}, we
review some conventions on Legendrian and transverse knots.  In
Section~\ref{sec:CombinatorialHFK}, we recall the combinatorial chain
complex from a grid diagram~\cite{MOS}. In
Section~\ref{sec:GridMoves}, we describe the relationship between grid
diagrams, Legendrian knots, and their stabilizations. In
Section~\ref{sec:CombinatorialMoves}, we recall the isomorphisms on
homology induced by basic moves on grid diagrams~\cite{MOST}.  The
Legendrian invariant is defined in
Section~\ref{sec:LegendrianElement}, and the properties stated above
are established there.
In this introduction, and indeed throughout most of this paper, we
have focused on the case of knots, rather than links, mainly for
notational simplicity. In Section~\ref{sec:Links}, we
extend the results to the case of links in $S^3$.
Finally, in Section~\ref{sec:examples} we give examples
of the computation of the Legendrian elements.  We reprove that the
knot~$5_2$ is not Legendrian simple, and show that the link $6^2_3$ is
not, either.

\subsection*{Acknowledgements}
The authors wish to thank John Baldwin, John
Etnyre, Ciprian Manolescu, Lenhard Ng, Jacob Rasmussen, and Andr{\'a}s
Stipsicz for helpful discussions during the course of this work. We
are especially indebted to Ng for sharing with us his expertise, and
clarifying some questions about symmetries of Legendrian knots.
Moreover, it appears that Baldwin and Ng have both independently
discovered at least some of the structures described here.


\section{Legendrian and transverse knots}
\label{sec:Contact}

We first recall some standard definitions from contact topology.  Etnyre's
survey~\cite{Etnyre} is a good reference for this material.

Endow $S^3$ with its standard contact structure $\xi$.  Restricted to
$\R^3\subset S^3$, this contact structure is the two-plane
distribution which is the kernel of the one-form $dz-y\,dx$.

Recall that a \emph{Legendrian knot} is a smooth knot $\LegK\subset
S^3$ whose
tangent vectors are contained in the contact planes of $\xi$. Two
knots are \emph{Legendrian isotopic} if they can be connected by a
smooth one-parameter family of Legendrian knots.

There are two classical invariants of Legendrian knots, the {\em
rotation number} and the \emph{Thurston\hyp Bennequin invariant}, as
follows.

Fix an embedded Seifert surface $F$
for the oriented Legendrian knot $\oLegK$. The restriction of $\xi$ to~$F$
determines an oriented two-plane bundle over~$F$, which has a
trivialization along the boundary induced by tangent vectors to the
knot. The rotation number~$\Rotation(\oLegK)$ is the relative first Chern
number of this two-plane field over $F$, relative to the
trivialization over $\partial F$.

On the other hand, the restriction of~$\xi$ to~$\LegK$ determines a
framing of~$\LegK$. The Thurston\hyp Bennequin number~$\TB(\LegK)$ is the
self-linking number of $\oLegK$ with respect to this framing.  That
is, if we let $\oLegK'$ denote a push-off of $\oLegK$ with respect to
the framing, the Thurston\hyp Bennequin invariant is the oriented
intersection number of $\oLegK'$ with $F$.

Note that the overall sign
of the rotation number depends on the choice of orientation for
$\oLegK$, but the Thurston\hyp Bennequin invariant is independent of this
choice.

It is customary to study Legendrian knots via their \emph{front
  projections}, defined by the projection map $(x,y,z)\mapsto (x,z)$.
The front projection of a Legendrian embedding has no vertical
tangencies, and in the generic case, its only singularities are
double-points and cusps.

The classical invariants of (oriented) Legendrian knots $\oLegK$ can
be read off from their (oriented) front projections $\Pi$; indeed, we have
\begin{align}
\TB(\LegK)&=\writhe(\Pi)-\OneHalf\#\{\text{cusps in $\Pi$}\}\label{eq:TBformula} \\
\Rotation(\oLegK)&=\OneHalf \Big(\#\{\text{downward-oriented cusps}\} -
\#\{\text{upward-oriented cusps}\}\Big) \label{eq:Rformula},
\end{align}
where here $\writhe(\Pi)$ denotes the \emph{writhe} of the projection,
i.e., the number of positive minus the number of
negative crossings of the projection.

A \emph{transverse knot} is a knot $\TransK\subset S^3$ is a knot whose
tangent vectors are transverse to the contact planes of $\xi$. Two
transverse knots are \emph{transverse isotopic} if they can be
connected by a smooth one-parameter family of transverse knots.

There is a classical invariant for transverse knots $\TransK$, the
\emph{self-linking number}, $\SL(\TransK)$. To define this, observe
first that since $\xi$ is cooriented, transverse knots inherit a
canonical orientation: the orientation for which the intersection
number of $\gamma'(t)$ with the plane $\xi_{\gamma(t)}$ is positive
for all $t\in[0,1]$ (or, equivalently, the evaluation of $dz-y\,dx$ is
positive on a vector oriented in the direction of the knot). Next,
consider a Seifert surface $F$ for $\TransK$ compatible with its
orientation. Now, along $\TransK$, the tangent space to $F$ and
the contact planes $\xi$ intersect in a line field, which in turn
inherits a natural orientation as ``outward pointing'' along
$F$. Thus, we have a trivialization of $\xi|_{\partial F}$. The self-linking
number, then, is the relative first Chern number of $\xi$ on $F$
relative to this trivialization of $\xi$ along its boundary.

\begin{figure}
\begin{center}
\graph{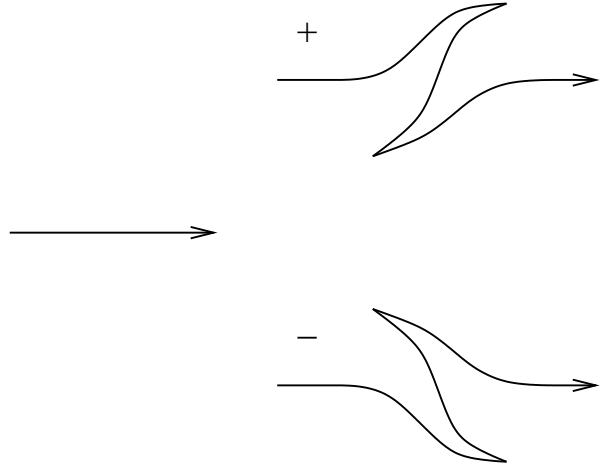}
\end{center}
\caption {{\bf Introducing cusps.}
\label{fig:StabilizeLegendrian}
        Given an oriented arc in an oriented Legendrian knot,
        we can introduce two cusps locally, in two ways, as pictured.
        The top is \emph{positive stabilization}, the bottom is
        \emph{negative stabilization}.}
\end{figure}

Given an oriented Legendrian knot~$\oLegK$, there are arbitrarily
close smooth curves (in the $C^\infty$ topology) which are transverse
(and which inherit the same orientation as~$\oLegK$).  We call such a
knot the (positive) \emph{transverse push-off} of~$\oLegK$.
Conversely, a transverse knot has a regular neighborhood which is
contactomorphic to the solid torus endowed with a standard contact
structure, for which the core is transverse.  There are Legendrian
curves in this solid torus which are transverse to the meridional
disks, meeting each disk in a single, transverse point of
intersection. We call these Legendrian curves the \emph{Legendrian
  approximations} to the transverse knot~$\TransK$. The transverse
push-off of a Legendrian approximation to~$\TransK$ is transversally
isotopic to~$\TransK$, and a Legendrian knot~$\oLegK$
becomes Legendrian isotopic, after it is negatively stabilized
sufficiently many times, to a Legendrian
approximation to its transverse push-off.  If $\TransK$ is the
transverse push-off of~$\oLegK$, then
\[
\SL(\TransK) = \TB(\oLegK) - \Rotation(\oLegK).
\]

Any two Legendrian approximations to a transverse knot can be
negatively stabilized so that they become Legendrian
isotopic~\cite{TransverseIsotopic}.  Thus, Legendrian knots modulo
negative stabilization are the same as transverse knots modulo
isotopy.


\section{Combinatorial knot Floer homology}
\label{sec:CombinatorialHFK}

We now recall the combinatorial chain complex for Heegaard-Floer
homology from earlier papers
\cite{MOS,MOST}, to which the reader is directed for a more in-depth
treatment.

Recall that a (planar) \emph{grid diagram} $G$ lies in a square grid on
the the plane with $n \times n$ squares. Each square is decorated
either with an $X$, an $O$, or nothing, arranged so that:
\begin{itemize}
\item every row contains exactly one $X$ and one $O$; and
\item every column contains exactly one $X$ and one $O$.
\end {itemize}
The number $n$ is called the \emph{grid number} of $G$.
We number the $O$'s and $X$'s by
$\{O_i\}_{i=1}^n$ and $\{X_i\}_{i=1}^n$, and we denote the
two sets by~$\Os$ and~$\Xs$, respectively.  (We use here the
notation from~\cite{MOST}; the $O_i$ correspond to the ``white dots''
of~\cite{MOS} and the $w_i$ of~\cite{Links}, while the $X_i$ to the
``black dots'' of~\cite{MOS} and the $z_i$ of~\cite{Links}.)

Given a planar grid diagram $G$, we place it in a standard
position on the plane by placing the bottom left corner at the
origin, and making each cell a square of edge length one. We then
construct an oriented, planar link projection by drawing horizontal
segments from the $O$'s to the $X$'s in each row, and vertical
segments from the $X$'s to the $O$'s in each column. At every
intersection point, we let the horizontal segment be the underpass and
the vertical one the overpass. This produces a planar diagram for an
oriented link $\orL$ in $S^3$. We say that $\orL$ has a grid
presentation given by $G$. It is easy to tell if a given grid
presentation determines a knot. For the moment, we restrict attention
to this case, and return to the more general case in
Section~\ref{sec:Links}.

If we cyclically permute the rows or column of a grid diagram, we
do not change the link that it represents, so from now on we will think
of grid diagrams as being drawn on a torus~$\Torus$.  Let the
horizontal (respectively
vertical) \emph{(grid) circles} be the circles in between two adjacent
rows (respectively columns) of marked squares.

Given a toroidal grid diagram~$G$, we associate to it a chain complex
$\bigl(\Cm(G;\Field), \dm \bigr)$ as follows.
Let $\S=\S(G)$ be the set of one-to-one correspondences between the
horizontal and vertical grid circles. More geometrically, we can think of
elements of $\S$ as $n$-tuples of intersection points between the
horizontal and vertical grid circles, with the property that there is
exactly one intersection point on each horizontal and vertical
grid circle.
Let $\Cm(G;\Field)$ be the free module over $\Field[U_1,\dots,U_n]$
generated by elements of $\S$, where here the $\{U_i\}_{i=1}^n$ are indeterminates.

The complex has a bigrading, induced by two functions $A \colon \S
\longrightarrow \Z$ and $M \colon \S \longrightarrow \Z$ defined as
follows.

Given two collections $A$, $B$ of finitely many points in the plane,
let $\NEoff(A,B)$ be the number of pairs $(a_1,a_2)\in A$ and
$(b_1,b_2)\in B$ with $a_1<b_1$ and $a_2<b_2$, and let $\NESW(A,B) =
(\NEoff(A,B)+\NEoff(B,A))/2$.  Take a fundamental
domain for the torus which is cut along a horizontal and vertical
circle, with the left and bottom edges included.  Given a generator
$\x\in\S$, we view $\x$ as a collection of points in this fundamental
domain.  Similarly, we view $\Os=\{O_i\}_{i=1}^n$ as a collection of
points in the plane.  Define
\[
M(\x):=\NESW(\x,\x)-2\NESW(\x,\Os) +\NESW(\Os,\Os)+1.
\]
We find it convenient to write this formula more succinctly as
\begin{equation}
  \label{eq:DefMaslov}
    M(\x)=\NESW(\x-\Os,\x-\Os)+1,
\end{equation}
where we extend $\NESW$ bilinearly over formal sums (or differences)
of subsets.
$M(\x)$ depends only on the sets $\x$ and $\Os$, but not on how we
drew the torus on the
plane, as we showed
earlier \cite[Lemma~\ref{OnComb:lemma:MaslovWellDefined}]{MOST}. Furthermore,
by the argument there we can alternately compute $M$ using a fundamental domain
that includes the right and top edges instead of the left and bottom edges.

Define $M_S(\x)$ to be the same as $M(\x)$ with the set~$S$ playing
the role of $\Os$.  We define
\begin{equation}\label{eq:AlexanderFormulaOne}
  \begin{split}
    A(\x) &:= \frac{1}{2}(M_{\Os}(\x) - M_{\Xs}(\x))
       - \Bigl(\frac{n-1}{2}\Bigr)\\
    &= \NESW(\x - \frac{1}{2}(\Xs + \Os), \Xs - \Os)
       - \Bigl(\frac{n-1}{2}\Bigr).
  \end{split}
\end{equation}

The module $\Cm(G;\Field)$ inherits a bigrading from the functions $M$ and $A$
above, with the additional convention that multiplication by $U_i$
drops the Maslov grading by two and the Alexander grading by one.

Given a pair of generators $\x$ and $\y$, and an embedded rectangle
$r$ in $\Torus$ whose edges are arcs in the horizontal and vertical
circles, we say that $r$ \emph{connects} $\x$ to $\y$ if $\x$ and $\y$ agree
along all but two horizontal circles, if all four corners of $r$ are
intersection points in $\x\cup\y$, and if when we traverse each horizontal
boundary component of $r$ in the direction dictated by the
orientation that $r$ inherits from $\Torus$, the arc is oriented
from a point in $\x$ to the point in $\y$.  Let $\Rect(\x,\y)$ denote
the collection of rectangles connecting $\x$ to $\y$.  If $\x,\y\in
\S$ agree along all but two horizontal circles, then there are exactly
two rectangles in $\Rect(\x,\y)$; otherwise $\Rect(\x,\y)=\emptyset$.
A rectangle~$r\in\Rect(\x,\y)$ is said to be \emph{empty} if
$\Interior(r)\cap\x = \emptyset$, or equivalently if
$\Interior(r)\cap\y = \emptyset$.  The space of empty rectangles
connecting~$\x$ and~$\y$ is denoted $\EmptyRect(\x,\y)$.

We endow $\Cm(G;\Field)$ with an endomorphism $\dm\colon
\Cm(G;\Field)\longrightarrow \Cm(G;\Field)$ defined by
\begin{equation}
\label{eq:DefCm}
\dm(\x)=\sum_{\y\in\S}\,
\sum_{r\in\EmptyRect(\x,\y)}\!\!
U_1^{O_1(r)}\cdots U_n^{O_n(r)}\cm \y,
\end{equation}
where
$O_i(r)$ denotes the number of times $O_i$ appears in the interior of~$r$.
This differential decreases the Maslov grading by~$1$ and preserves a
filtration by the Alexander grading.

The filtered chain homotopy type of the above
complex is a knot
invariant~\cite{MOS}; indeed,
it is the filtered Heegaard-Floer complex~\cite{Knots,RasmussenThesis}. In
the sequel, we will need a slightly less refined version of this. Let
$\CKm(G;\Field)$ be the associated graded object. That is,  $\CKm(G;\Field)$ has the
same bigraded set of generators as $\Cm(G;\Field)$, but its differential now
counts only empty rectangles with no elements of $\Xs$ in them; i.e., it
is given by 
\begin{equation}
\label{eq:DefCLm}
\partial(\x)=\sum_{\y\in\S}\,
\sum_{\substack{r\in\EmptyRect(\x,\y)\\r \cap \Xs = \emptyset}}\!\!
    U_1^{O_1(r)}\cdots U_n^{O_n(r)}\cm \y.
\end{equation}
This complex is now bigraded, splitting as
$\CKm(K;\Field)=\bigoplus_{m,a\in\Z}\CKm_m(K,a;\Field)$, so that the
differential
drops the Maslov grading~$m$ by one and preserves the Alexander
grading~$a$. Thus, its homology groups inherit a bigrading
$$\HFKm(K;\Field)=\bigoplus_{m,a\in\Z}\HFKm_m(K,a;\Field).$$
These groups can be viewed as a bigraded module over $\Field[U]$,
where $U$ acts by multiplication by $U_i$ for any $i=1,\dots,n$.

\begin{theorem}[\cite{MOS}]
        \label{thm:MOS}
  The bigraded homology groups $\HFKm(K;\Field)$ agree with the knot Floer
  homology of~$K$ with coefficients in~$\Field$.
\end{theorem}

Moreover, the chain complex can be lifted to $\Z$
coefficients.  Let $\Cm(G)$ be the free $\Z[U_1,\ldots,U_n]$ module
generated by~$\S$.  Then we have:

\begin{theorem}[{\cite[Theorem~\ref{OnComb:thm:SignAssignments}]{MOST}}]
\label{thm:MOST} There is an essentially unique function $\sign\colon
\Rect(\x,\y)\rightarrow \{\pm 1\}$ with the following two
properties:
\begin{itemize}
\item the endomorphism $\dm\colon \Cm(G) \longrightarrow \Cm(G)$
defined by
$$\dm_\sign(\x)=\sum_{\y\in\S}\,
\sum_{r\in\EmptyRect(\x,\y)}\!\!
\sign(r) \cm U_1^{O_1(r)}\cdots U_n^{O_n(r)}\cm \y
$$
is a differential, and
\item the homology of the associated complex when we
set all $U_i=1$ has non-zero rank; i.e.,
$H_*(\Cm(G)/\{U_i=1\}_{i=1}^n)\otimes {\mathbb Q}$ is non-zero.
\end{itemize}
Moreover, the filtered quasi-isomorphism type of the complex
$(\Cm(G),\dm_\sign)$, thought of as a complex over $\Z[U]$ (where $U$
acts as multiplication by any $U_i$),
is an invariant of the link.  In particular, it is independent of the
choice of $\sign$ with the above properties and the grid diagram for~$K$.
\end{theorem}

The proof of Theorem~\ref{thm:MOST} gives an independent (and
elementary) proof that
the homology groups are a topological invariant of~$K$.

\begin{remark}
The group $\HFKm(K)$ is the homology group of the associated graded
object of the filtered complex $\Cm(K)$. The filtered
quasi-isomorphism type of~$\Cm(K)$ is a more refined knot
invariant, and indeed the more general version of
Theorem~\ref{thm:MOS} identifies this filtered quasi-isomorphism type
with a corresponding more general object associated to knots defined
using holomorphic disks. In more concrete terms, as a result of this
extra structure, $\HFKm(K)$ is endowed with a collection of higher
differentials, the first of which is
$$\delta_1 \colon \HFKm_{d}(K,s) \longrightarrow \HFKm_{d-1}(K,s-1).$$
The identification of, say, Theorem~\ref{thm:LegendrianInvariant}
induce isomorphisms
$$\phi\colon \HFKm_{d}(K(G),s) \longrightarrow \HFKm_{d}(K(G'),s)$$
which commute with $\delta_1$. We have no need for this extra structure
in the current paper, but in later work with Ng we use it to
distinguish distinct transverse knots \cite[Section 3.2]{NOT}.
(We encounter the filtered chain homotopy type briefly
in Section~\ref{sec:OnTau}, when comparing the usual definition of $\tau$
with the one discussed in the introduction.)
\end{remark}


\section{Grid diagrams and Legendrian knots}
\label{sec:GridMoves}

Grid presentations~$G$ can represent Legendrian or transverse
knots in addition to ordinary knots. Specifically, given a grid
presentation~$G$ of~$K$, we can construct a front projection for a
Legendrian realization of the mirror~$\Mirror(K)$ of~$K$ as follows.
Consider the projection of $K$ obtained from $G$ as in the previous section. It
is a projection with corner points, and indeed there are four types of
corner points, which we denote northwest, southwest, southeast and
northeast. Smooth all the northwest and southeast corners of the
projection, view the southwest and northeast corners as cusps, and
then tilt the diagram $45^\circ$ clockwise, so that the NE
(respectively SW)
corners become right (respectively left) cusps.
This gives a Legendrian front
projection for the mirror of the knot~$K$ described by~$G$. (It is
easy to find a different convention which does not give a mirror; the
present conventions appear to fit neatly with conventions on the
contact element~\cite{HolDiskContact}.)

Before giving the combinatorial presentation of Legendrian or
transverse knots, we first recall the combinatorial presentation of
links using grid diagrams.  There are several elementary moves on a
grid diagram~$G$ that do not change the topological link type:
\begin{description}
\item[Cyclic permutation]
  Cyclically permute the rows or columns of~$G$.
\item[Commutation]  For any pair of consecutive columns of~$G$
  so that the~$X$ and~$O$ from one column do not separate the~$X$
  and~$O$ on the other column, switch the decorations of these two
  columns, as in Figure~\ref{fig:CommutePicture}.
  In particular, if the $O$ and $X$ in one column are in adjacent
  rows, this move can be applied unless there is an~$X$ or an~$O$ in
  one of the same rows in the adjacent column.
  There is also a
  similar move where the roles of columns and rows are interchanged.
\item[Destabilization] For a corner~$c$ which is shared by a pair of
  vertically-stacked squares~$X_1$, $O_1$ marked with an $X$ and $O$
  respectively, we remove the markings of $X_1$ and $O_1$ and delete the
  horizontal and vertical circles containing~$c$. Indeed, we can (and
  do) assume that one of~$X_1$ or~$O_1$ meets an additional square
  marked by an $O$ or an~$X$, by the comment above about
  commutation when there is an adjacent $X$ and $O$. We also can (and
  do) assume that
  $c$ is either the lower-left or upper-right corner of~$O_1$.
\item[Stabilization] The inverse of destabilization.
\end{description}

\begin{prop}[Cromwell~\cite{Cromwell}, see also Dynnikov~\cite{Dynnikov}]
  Two grid diagrams represent the same topological link if and only if
  they can be connected by a sequence of cyclic permutation,
  commutation, stabilization, and destabilization moves.
\end{prop}

\begin{figure}
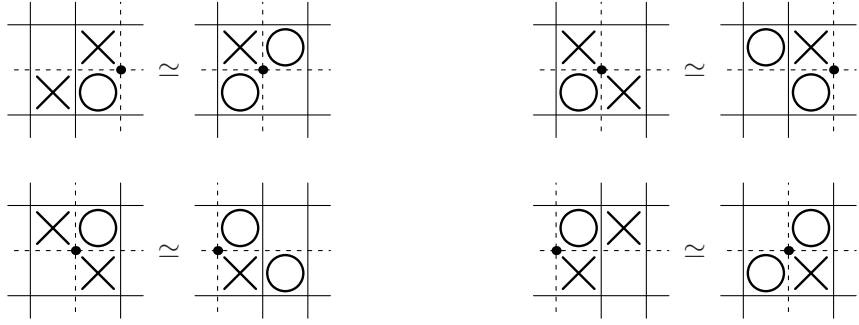

\begin{align*}
  \mfigb{destabs.0}&\simeq\mfigb{destabs.13} &
     \mfigb{destabs.1}&\simeq\mfigb{destabs.12} \\[8pt]
  \mfigb{destabs.2}&\simeq\mfigb{destabs.11} &
     \mfigb{destabs.3}&\simeq\mfigb{destabs.10}
\end{align*}
\caption {{\bf Destabilizations.}
\label{fig:Stabilizations}
We have enumerated here the eight types of destabilizations.
The dotted lines are to be removed in the destabilized picture.
Starting from the upper left corner and reading along the top row first,
we have destabilizations of types:
$\XNW$, $\OSE$, $\XNE$, $\OSW$,
$\XSW$, $\ONE$, $\XSE$, and $\ONW$. In each destabilization the two dotted
lines are removed and the corner~$c$ is marked. The indicated pairs of
destabilizations are equivalent modulo commutation moves on the torus.}
\end{figure}

We can further classify (de)stabilization moves according to the local
configuration of $X$'s and $O$'s. Recall that we are assuming now that 
three marked squares in the original diagram share one corner. 
There are two data to keep track
of: the marking shared by two of these three squares (i.e., an $X$ or an $O$),
and the 
placement of the \emph{unmarked} square relative to the shared corner,
either $\NW$,
$\SW$, $\SE$, or $\NE$.
See Figure~\ref{fig:Stabilizations}. 

The following two lemmas follow by elementary manipulations.

\begin{lemma}
  \label{lem:stab-types}
  A stabilization of type~$\OSE$ (respectively $\ONE$, $\ONW$, or $\OSW$) is
  equivalent to a stabilization of type~$\XNW$ (respectively $\XSW$, $\XSE$,
  or $\XNE$) followed by a sequence of commutation moves on the torus.
\end{lemma}
\begin{proof}
  After a stabilization at an $X$ vertex, we can slide either of the
  resulting segments of length~1 to a neighboring vertex (of
  type~$O$) by a sequence of commutation moves.  A straightforward
  check shows that we get a stabilization of the type as indicated in
  the statement.
\end{proof}

\begin{lemma}
  \label{lem:cyclic-perm}
  A cyclic permutation is equivalent to a sequence of commutations in
  the plane and
  (de)stab\-il\-i\-za\-tions of types~$\XNW$, $\XSE$, $\ONW$, and $\OSE$.
\end{lemma}

\begin{proof}
  Since the allowed moves are symmetric under reversing orientation,
  we may suppose without loss of generality that we wish to move a
  horizontal segment from
  the top to the bottom, with left end of the segment marked~$X_1$ and right
  end marked~$O_2$.  Let~$O_1$ (respectively $X_2$) be the other mark in the
  column containing~$X_1$ (respectively $O_2$).  Apply a
  stabilization of type $\XNW$ at~$X_2$, and commute
  the resulting horizontal segment of length~1 to the bottom of the
  diagram.  We now have a vertical segment stretching the height of
  the diagram; apply commutation moves until it is just to the left of
  the column containing~$X_1$.  Now the horizontal segment starting
  at~$X_1$ is of length~1, and so can be commuted down until it is
  just above~$O_1$, where we can apply a destabilization of type
  $\OSE$ to get the desired cyclic permutation.
\end{proof}

\begin{prop}
  \label{prop:LegendrianGrid}
  Two grid diagrams represent the same Legendrian link if and only if
  they can be connected by a sequence of commutation and
  (de)stabilizations of types~$\XNW$ and~$\XSE$ on the torus.
\end{prop}

\begin{proof}
  By Lemmas~\ref{lem:stab-types} and~\ref{lem:cyclic-perm},
  we can  equivalently consider commutation and
  (de)stab\-il\-i\-za\-tions of all types $\NW$ and $\SE$ in the
  rectangle (rather than on the torus).

  We must now check that each type of commutation and allowed
  stabilization (in the rectangle) gives an isotopy of the
  corresponding Legendrian knot.  Indeed, after rotating $45\deg$ and
  turning the corners into smooth turns or cusps as appropriate, each
  elementary move of the grid diagram becomes a sequence of Legendrian
  Reidemeister moves of the front projection. For instance, as shown in
  Figure~\ref{fig:nw-stab}, an $\SE$ stabilization becomes either a
  planar isotopy or a Legendrian Reidemeister~1 move, depending on the
  relation of the stabilized corner to the rest of the diagram.  An
  example of a commutation move that turns into a Legendrian
  Reidemeister~2 move is shown in Figure~\ref{fig:type-exchange-move};
  other commutation moves are similar, and may also involve
  Reidemeister~3 moves.

  \begin{figure}
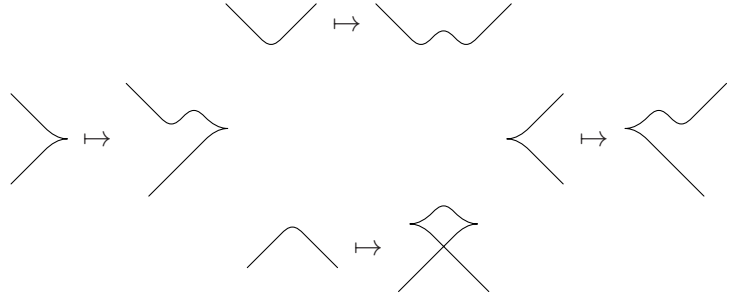

    $\mfigb{legendrian.11}\mapsto\mfigb{legendrian.21}$\\[8pt]
    $\mfigb{legendrian.10}\mapsto\mfigb{legendrian.20}\hspace{3.5cm}
      \mfigb{legendrian.12}\mapsto\mfigb{legendrian.22}$\\[-3pt]
    $\mfigb{legendrian.13}\mapsto\mfigb{legendrian.23}$\\
    \caption{The four different ways an $\SE$ stabilization can appear
      after converting into a Legendrian front.
      Three of them are planar isotopy, while the fourth is an allowed
      Reidemeister move on Legendrian knots.}
    \label{fig:nw-stab}
  \end{figure}

  \begin{figure}
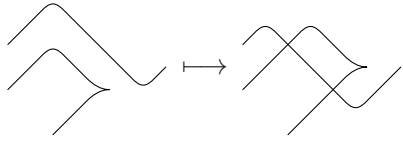

    \[
    \mfigb{legendrian.30} \longmapsto \mfigb{legendrian.31}
    \]
    \caption{A commutation move giving a Legendrian Reidemeister~2 move.}
    \label{fig:type-exchange-move}
  \end{figure}

  To go the other direction, we must show that every sequence of
  Reidemeister moves on a Legendrian front can be turned into a
  corresponding sequence of grid moves.  First note that any
  Legendrian front projection can be turned into a grid diagram: Take
  the Legendrian front and stretch it horizontally until no portion of
  the diagram is at an angle of more than $45\deg$ from the
  horizontal.  Then the curve can be approximated by a sequence of
  straight segments at an angle of $\pm 45\deg$.  After rotating by
  $45\deg$ counter\hyp clockwise and adjusting the segments to have
  consecutive integer coordinates, we have a grid diagram
  corresponding to the front projection.

  We can do the same thing with any Legendrian isotopy: stretch the
  intermediate diagrams so that no edges are too steep and approximate
  each one by a sequence of straight segments.  It is an elementary
  verification that each modification along the way (i.e., change of the
  approximation by segments and Legendrian Reidemeister moves) can be
  achieved by
  a sequence of commutation and allowed destabilizations.
\end{proof}

Furthermore, a stabilization of type~$\XSW$ is
a negative stabilization of the Legendrian link (see
Figure~\ref{fig:stab-xsw}), so we have:

\begin{figure}
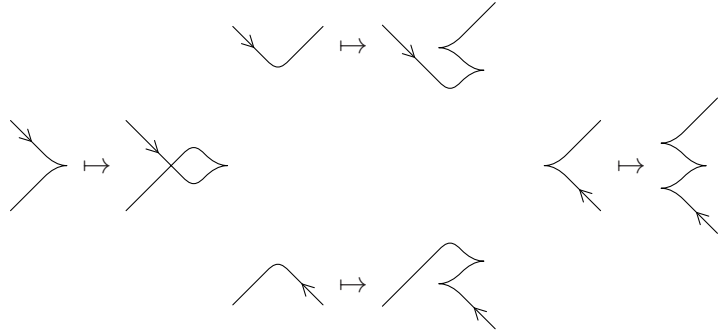

  $\mfigb{legendrian.16}\mapsto\mfigb{legendrian.26}$\\[-3pt]
  $\mfigb{legendrian.15}\mapsto\mfigb{legendrian.25}\hspace{4cm}
    \mfigb{legendrian.17}\mapsto\mfigb{legendrian.27}$\\[-3pt]
  $\mfigb{legendrian.18}\mapsto\mfigb{legendrian.28}$\\
  \caption{Stabilizations of type $\XSW$, after conversion to a Legendrian front.  In each case the Legendrian knot type is changed by a negative stabilization.}
  \label{fig:stab-xsw}
\end{figure}

\begin{corollary}
  Two grid diagrams represent the same transverse link if and only if
  they can be connected by a sequence of commutation and
  (de)stabilizations of types $\XNW$, $\XSE$, and $\XSW$.
\end{corollary}

Finally we note how symmetries of the grid diagram relate to
symmetries of the knot.  We consider only the symmetries that
preserve the set of cusps.

\begin{lemma}
  \label{lem:GridSymmetries}
  Symmetries of a grid diagram~$G$ take an oriented Legendrian
  knot~$\oLegK$ to the following Legendrian knots:
  \begin{itemize}
  \item Reflection through the $x=-y$ axis: $-\oLegK$;
  \item Reflection through the $x=y$ axis: $-\mu(\oLegK)$; and
  \item Rotation by $180\deg$: $\mu(\oLegK)$.
  \end{itemize}
\end{lemma}

\begin{proof}
  Reflection of~$G$ through the $x=-y$ axis corresponds to reflecting
  the front projection through the $z$ (vertical) axis, which in turn
  corresponds to rotating the Legendrian knot by $180\deg$ around the
  $z$-axis, which is a Legendrian isotopy.  (The fact this is a
  Legendrian isotopy is most easily seen
  by using the isotopic contact form $dz - (x\,dy-y\,dx)/2$, which is
  rotationally symmetric.)  Moreover, by considering the orientation
  conventions for knot diagrams, one can easily see that this
  operation also reverses the orientation of the knot.
  
  Reflection of~$G$ through the $x=y$ axis similarly corresponds to
  rotation of the Legendrian projection by~$180\deg$ around the axis
  normal to the Legendrian projection (which is the definition of
  Legendrian mirror), followed by a reflection through the $x=-y$ axis
  (which does not affect the Legendrian isotopy class). This operation
  also reverses the orientation.

  Rotation of~$G$ by $180\deg$ is the composition of the previous two
  symmetries.
\end{proof}


\section{Grid moves for knot Floer homology}
\label{sec:CombinatorialMoves}

We recall here the explicit maps induced by commutation and
destabilization moves in the combinatorial proof of topological
invariance of $\HFKm$~\cite{MOST}.  We will give the formulas
over~$\Field$; signs can be added to all these formulas to give
formulas that work
over~$\Z$~\cite[Section~\ref{OnComb:subsec:DiscussionOverZ}]{MOST}.

\subsection{Commutation maps}
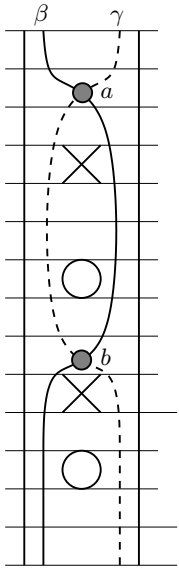
\begin{figure}
\begin{center}
\begin{picture}(0,0)%
\graph{CommutePicture.pstex}%
\end{picture}%
\setlength{\unitlength}{1579sp}%
\begingroup\makeatletter\ifx\SetFigFont\undefined%
\gdef\SetFigFont#1#2#3#4#5{%
  \reset@font\fontsize{#1}{#2pt}%
  \fontfamily{#3}\fontseries{#4}\fontshape{#5}%
  \selectfont}%
\fi\endgroup%
\begin{picture}(2724,8886)(1189,-7294)
\put(2701,-4111){\makebox(0,0)[lb]{\smash{{\SetFigFont{9}{10.8}{\rmdefault}{\mddefault}{\updefault}{\color[rgb]{0,0,0}$b$}%
}}}}
\put(2701, 89){\makebox(0,0)[lb]{\smash{{\SetFigFont{9}{10.8}{\rmdefault}{\mddefault}{\updefault}{\color[rgb]{0,0,0}$a$}%
}}}}
\put(1651,1289){\makebox(0,0)[lb]{\smash{{\SetFigFont{9}{10.8}{\rmdefault}{\mddefault}{\updefault}{\color[rgb]{0,0,0}$\beta$}%
}}}}
\put(2851,1289){\makebox(0,0)[lb]{\smash{{\SetFigFont{9}{10.8}{\rmdefault}{\mddefault}{\updefault}{\color[rgb]{0,0,0}$\gamma$}%
}}}}
\end{picture}%
\end{center}
\caption {{\bf Commutation.}
A commutation move, viewed as replacing one vertical circle ($\beta$, 
undashed) with another ($\gamma$, dashed).}
\label{fig:CommutePicture}
\end{figure}

More explicitly, suppose that $G$ and $H$ are two grid diagrams for
the same oriented knot $\oKnot$, which differ by commuting two vertical
edges. It is convenient to draw both diagrams on the same torus,
replacing a distinguished vertical circle $\beta$ for $G$ with a
different one $\gamma$ for $H$, as pictured in
Figure~\ref{fig:CommutePicture}.  The circles $\beta$ and $\gamma$
meet each other transversally in two points~$a$ and~$b$, which are not
on a horizontal circle.

We define a chain map $\Phi_{\beta\gamma} \colon \Cm(G)
\longrightarrow \Cm(H)$ by counting pentagons in the torus.  Given
$\x\in \S(G)$ and $\y\in\S(H)$, we let $\Pent_{\beta\gamma}(\x,\y)$
denote the space of embedded pentagons with the following properties,
as illustrated in Figure~\ref{fig:FindPentagons}.
This space is empty unless $\x$ and $\y$ coincide at $n-2$ points. An
element of $\Pent_{\beta\gamma}(\x,\y)$ is an embedded disk in
$\Torus$, whose boundary consists of five arcs, each contained in
horizontal or vertical circles. Moreover, under the orientation
induced on the boundary of $p$, we start at the $\beta$-component of
$\x$, traverse the arc of a horizontal circle, meet its corresponding
component of $\y$, proceed along an arc of a vertical circle, meet the
corresponding component of $\x$, continue through another
horizontal circle, meet the component of $\y$ contained in the
distinguished circle $\gamma$, proceed along an arc in $\gamma$, meet an
intersection point of $\beta$ with $\gamma$, and finally, traverse an
arc in $\beta$ until we arrive back at the initial component of $\x$.
Finally, all the angles here are required to be acute. These
conditions imply that there is a particular intersection point,
denoted $a$, between $\beta$ and $\gamma$ which appears as one of the
corners of any pentagon in $\Pent_{\beta\gamma}(\x,\y)$.  The other
intersection point $b$ appears in all of the pentagons in
$\Pent_{\gamma\beta}(\y,\x)$.  The space of empty pentagons $p\in
\Pent_{\beta\gamma}(\x,\y)$ with $\x \cap \Interior(p) = \emptyset$,
is denoted $\EmptyPent_{\beta\gamma}$.

Given $\x\in\S(G)$, define
\begin{align*}
\Phi_{\beta\gamma}(\x) =& \sum_{\y\in \S(H)}\,\,
    \sum_{p\in\EmptyPent_{\beta\gamma}(\x,\y)}\!\!
U_1^{O_1(p)}\cdots U_n^{O_n(p)} \cm \y
\in \Cm(H).
\end{align*}

It is elementary to see that the above map induces a chain homotopy
equivalence~\cite[Proposition~\ref{OnComb:prop:Commute}]{MOST}.

In this paper we will consider the above map on the associated graded
object $\CKm(K)$, i.e., where we count $p\in \EmptyPent_{\beta\gamma}$
subject to the further constraint that $p \cap \Xs = \emptyset$.

\subsection{Stabilization maps}
\label{subsec:Stabilization}
Next, we consider the stabilization map.  Let $G$ be a grid diagram
and $H$ denote a stabilization.  We discuss in detail the case where
we introduce a new column with $O_1$ immediately above $X_1$, and
there is another marking~$X_2$ immediately to the left or to the right
of $O_1$, as is the
case of two of the four types of $X$\hyp stabilization; the cases where
$X_1$ is immediately above $O_1$ can be treated symmetrically by a
rotation of all diagrams by $180^\circ$.

Label the $O$ in the same row as $X_1$ by $O_2$. Let $\beta_1$ be the vertical
circle just to the left of $O_1$ and $X_1$, and let $\alpha$ denote the horizontal
circle separating the squares marked~$O_1$ and~$X_1$.

Let $B=\Cm(G)$, $C=\Cm(H)$, and let $C'$ be the mapping cone of
$$U_2-U_1\colon B[U_1]\longrightarrow B[U_1],$$
i.e.,
$C'[U_1]=B[U_1]\oplus B[U_1]$, endowed with the differential
$\partial'\colon C' \longrightarrow C'$ given by
$$
\partial'(a,b)=(\dm a, (U_2-U_1) \cm a-\dm b) $$
where
here $\dm$ denotes the differential within $B$.  Note that
$B$ is a chain complex over $\Z[U_2,\ldots,U_n]$, so that $B[U_1]$
denotes the induced complex over $\Z[U_1,\ldots,U_n]$ gotten by
introducing a new formal variable $U_1$.  Let $\Left\cong B[U_1]$
(respectively
$\Right$) be the subgroup of $C'$ of elements of the form
$(c,0)$ (respectively $(0,c)$) for $c\in B[U_1]$.
The module $\Right$ inherits
Alexander and Maslov gradings from its identification with $B[U_1]$,
while $\Left$ is given the Alexander and Maslov gradings which are one
less than those it inherits from its identification with $B[U_1]$.
With respect to these conventions, the mapping cone is a filtered
complex of $\Ring$-modules. 
We claim that $C'$ is quasi-isomorphic to $B$. In fact, it is
straightforward to verify that a quasi-isomorphism is given by the map
$(a,b)\mapsto b'$, where here $b'$ denotes the element of $B$ gotten
by taking $b\in B[U_1]$, and substituting $U_2$ for the formal
variable $U_1$.

Furthermore, there is a filtered quasi-isomorphism
\begin{equation}
  \label{eq:StabilizationMap}
  F\colon C \longrightarrow C'.
\end{equation}
To describe this, we introduce a little more notation.

Let $x_0$ be the intersection point of $\alpha$ and $\beta_1$. Let
$\Interval\subset \S(H)$ be the set of
$\x\in\S(H)$ which contain $x_0.$ There is, of course, a natural
(point-wise) identification between $\S(G)$ and $\Interval$, which
drops Alexander and Maslov grading by one.

As such, the differentials within $\Left$ and $\Right$ count
rectangles in~$H$
which do not contain $x_0$ on their boundary, although they may
contain $x_0$ in their interior.  Note however that the boundary
operator for rectangles containing $x_0$ does not involve the
variable~$U_1$.

\begin{figure}
\begin{center}
\graph{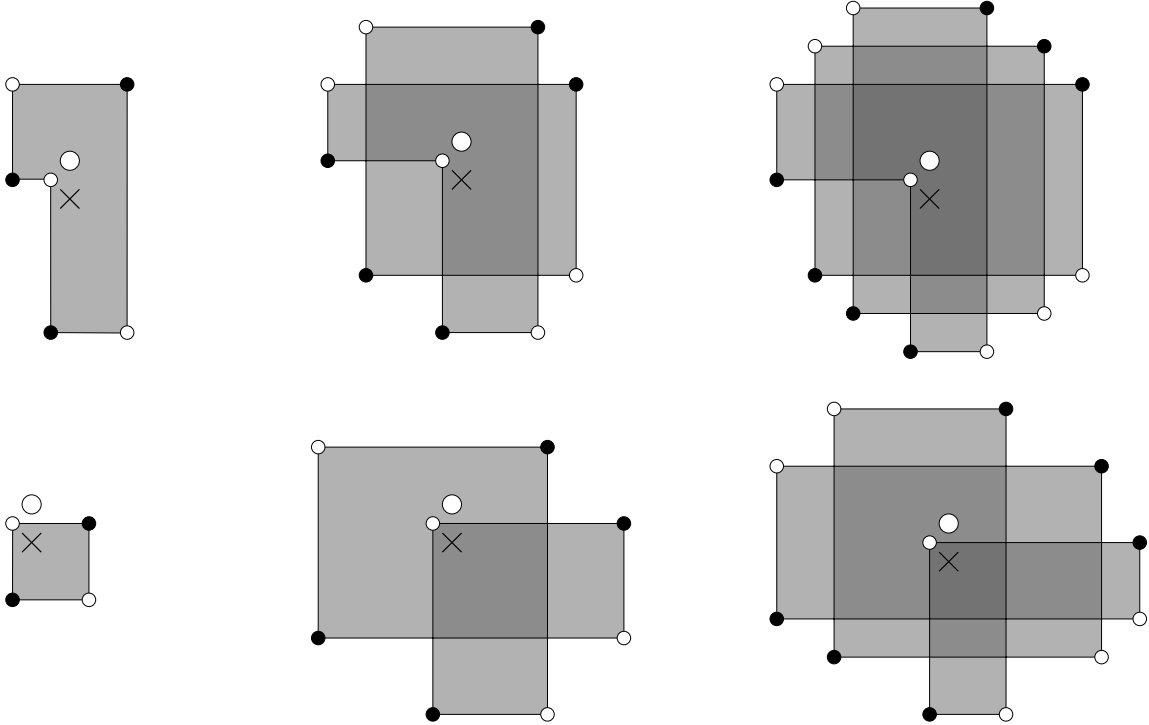}
\end{center}
\caption {{\bf Types of domains.}
We have listed here several domains in the stabilized diagram, labeling the
initial points by dark circles, and terminal points by empty circles.
The top row lists domains of type $L$, while the second row lists
some of type $R$. The marked $O$ and $X$ are the new ones in the
stabilized picture.
Darker shading corresponds to
higher local multiplicities.}
\label{fig:Domains}
\end{figure}

For the definition of $F$, we must consider objects more general than
rectangles, called domains. To define them, let us view the torus
$\Torus$ as a two-dimensional cell complex, with the toroidal grid
diagram inducing the cell decomposition with $n^2$ zero-cells, $2n^2$
one-cells and $n^2$ two-cells (the little squares). Let $U_{\alpha}$
be the one-dimensional sub-complex of $\Torus$ consisting of the union
of the $n$ horizontal circles.

\begin {definition}
Given $\x,\y\in \S(H)$, a \textit{path} from~$\x$ to~$\y$ is a
$1$--cycle $\gamma$ on the cell complex~$\Torus$, such that
the boundary of the intersection of $\gamma$ with $U_{\alpha}$ is
$\y-\x$.  
A \textit{domain}~$p$ from $\x$ to $\y$ is a two-chain in $\Torus$
whose boundary $\partial p$ is a path from $\x$ to $\y$.
The \emph{multiplicity} of $O_i$ in a domain~$p$, denoted $O_i(p)$, is
the the local multiplicity of the chain~$p$ at~$O_i$.
\end {definition}

For example, if $p$ is a rectangle from $\x$ to $\y$, then the above
definition of $O_i(r)$ coincides with the earlier definition used in
Equation~\eqref{eq:DefCm}.

Only certain domains will be counted in the definition of $F$.

\begin{definition}
  For $\x \in \S(H)$ and $\y \in \Interval \subset \S(H),$ a domain
  $p\in \pi(\x,\y)$ is said to be of type $L$ (respectively $R$) if either
  it is the zero chain, in which case $p$ has type~$L$, or it satisfies the
  following conditions:
  \begin{itemize}
  \item $p$ has only non-negative local multiplicities.
  \item For each $c\in\x\cup\y$, other than $x_0$, at least three of
    the four adjoining squares have vanishing local multiplicities.
  \item In a neighborhood of $x_0$ the local multiplicities in three
    of the adjoining rectangles are the same number~$k$.  When $p$ has
    type $L$, the lower left corner has local multiplicity~$k-1$,
    while for $p$ of type $R$ the lower right corner has
    multiplicity~$k+1$.
  \item $\partial p$ is connected.
  \end{itemize}
  The set of type $L$ (respectively $R$) domains from $\x$ to $\y$ is
  denoted $\pL(\x,\y)$ (respectively $\pR(\x,\y)$).  See
  Figure~\ref{fig:Domains} for examples.
\end{definition}

We now define maps
\begin{align*}
  \FL&\colon C\longrightarrow \Left\\
  \FR&\colon C\longrightarrow \Right
\end{align*}
where $\FL$ (respectively $\FR$) counts domains of type $L$
(respectively $R$)
without factors of $U_1$.  Specifically, define
\begin{align*}
\FL(\x)&=\sum_{\y\in\Interval}\,\sum_{p\in\pL(\x,\y)}
U_2^{O_2(p)}\cdots U_n^{O_n(p)}\cm\y \\
\FR(\x)&=\sum_{\y\in\Interval}\,\sum_{p\in\pR(\x,\y)}
U_2^{O_2(p)}\cdots U_n^{O_n(p)}\cm\y.
\end{align*}
In order to identify the range of $\FL$ (respectively $\FR$) with
$\Left$ (respectively $\Right$),
we implicitly use the identification $\S(G)\cong\Interval
\subset \S(H)$.

We put these together to define a map
\[
F=
\begin{pmatrix}
\FL \\ \FR 
\end{pmatrix} \colon C \longrightarrow C'.
\]
The fact that $F$ is a quasi-isomorphism is established
in~\cite[Proposition~\ref{OnComb:prop:Stabilization}]{MOST}.

Again, $F$ induces also a quasi-isomorphism on the associated graded
object, giving a map from $\CKm(H)$ to the mapping cone of
multiplication by $U_1-U_2$, thought of as an endomorphism of
$\CKm(G)[U_1]$. This induced map counts only those domains $p$ of type
$F$ for which $X_i(p)=0$ for all $i=2,\dots,n$.  (Note that we do allow
$X_1(p) \ne 0$.)


\section{Definition and invariance properties of the Legendrian invariants}
\label{sec:LegendrianElement}

We this setup, we can now construct the Legendrian invariant for
knots. The case of links works with minor modifications, as spelled out in 
Section~\ref{sec:Links}.

\begin{definition}
  \label{def:Canonical} Let $G$ be a grid diagram for a knot, and
  consider the chain complex $\CKm(G)$.  Consider elements $\xUR$,
  $\xLL\in \S(G)$ defined as follows.  Each component of $\xUR$ is the
  upper right corner of some square decorated with $X$, while each
  component of $\xLL$ is the lower left corner of some square
  decorated with $X$.  The chains $\zUR$ and $\zLL$ are defined to be
  $\pm\xUR$ and $\pm\xLL$, respectively, with signs specified in
  Definition~\ref{def:CanonicalWithSigns}.
\end{definition}

We will defer most discussion of signs until later.  For now, all the
proofs will work with an arbitrary choice of signs, although at
present we will
only prove some of the results up to a choice of sign. We now verify that
the above chains are in fact cycles.

\begin{lemma}
        \label{lemma:IsACycle}
        The elements $\zUR$ and $\zLL$ are  cycles in the chain
        complex $\CKm(G)$.
\end{lemma}

\begin{proof}
Consider any $y \in \S(G)$ and $r\in\Rect(\xUR,\y)$. Let $x_1\in\xUR$
be the upper right
corner of~$r$. By the definition of $\xUR$, there is an $X$ in the
square to the lower left of~$x_1$, and hence also an~$X$
in~$r$. Thus, $r$ cannot count in the definition of the
differential~$\partial$ for the associated graded object.
An analogous argument applies to~$\xLL$.
\end{proof}

Next, we calculate the Maslov and Alexander gradings.

\begin{lemma}
\label{lemma:AlexanderGrading}
        We have that $A(\xUR)=\frac{1}{2}M(\xUR)$
        and $A(\xLL)=\frac{1}{2}M(\xLL)$.
\end{lemma}

\begin{proof}
  Recall that
  $A(\xUR)=\frac{1}{2}(M_{\Os}(\xUR)-M_{\Xs}(\xUR))-\left(\frac{n-1}{2}\right)$.
  Because of the close relationship of~$\xUR$ and~$\Xs$, we have some
  equalities on the terms appearing in~$M_{\Xs}$ when we work in a
  fundamental domain including the right and top edges:
  \begin{align*}
    \NEoff(\xUR,\xUR)&=\NEoff(\Xs,\Xs)=\NEoff(\xUR,\Xs) \\
    \NEoff(\Xs,\xUR)&=\NEoff(\Xs,\Xs)+n.
  \end{align*}
  Therefore $M_{\Xs}(\xUR) = -n+1$ and $A(\xUR) = \frac{1}{2}M_{\Os}(\xUR)$ as
  desired. A similar argument applies to~$\xLL$.
\end{proof}

Let $K$ be the knot diagram represented by the grid diagram.  Remember
that the topological type of $K$ is the mirror of the topological type
of $\oLegK$.
\begin{lemma}
        \label{lemma:Maslov}
        The Maslov gradings of the elements $\xUR$ and $\xLL$ are given by 
        \begin{align}
          M(\xUR)&=-\writhe(K)-\#\{\text{downward-oriented cusps}\}+1 \label{eq:MasUR} \\
          M(\xLL)&=-\writhe(K)-\#\{\text{upward-oriented cusps}\}+1.  \label{eq:MasLL}
       \end{align}
\end{lemma}
        
\begin{proof}
  Each horizontal segment~$K_i$ of~$K$ goes from some $O_i$ to
  some $X_i$. Let $x_i$ be the point in $\xUR$ to the upper right
  of~$X_i$.  We claim that the quantity $C_i$ defined by
        $$C_i:=\NESW(\{x_i\} - \{O_i\},\xUR - \Os)$$
        is given by
        \begin{equation}
        \label{eq:GradingPerArc}
        \begin{split}
        C_i&=\#\{\text{negative crossings on $K_i$}\}
          -\#\{\text{positive crossings on $K_i$}\} \\
          &\qquad-\#\{\text{downward-oriented cusps among $\{X_i, O_i\}$}\}.
        \end{split}
        \end{equation}
        
        To prove this, we use the horizontal segment to divide up the
        plane into four
        regions, as follows.  Let $A$ be the vertical column
        through~$X_i$, $B$~be the vertical column through~$O_i$,
        $C$~be the vertical strip between~$A$ and~$B$ (note that the
        region $C$ can be empty when the arc from $O_i$ to $X_i$ has
        length one), and $D$ be the complement of $A\cup B\cup C$.
        (Note that $D$ generally has two connected components.)  Each
        region
        should be interpreted as including its right boundary.  Then
        for each $O_j$ in $D$,
        $\NESW(\{x_i\},\{O_j\})=\NESW(\{O_i\},\{O_j\})$, and similarly
        for each $x_j$ in $D$,
        $\NESW(\{x_i\},\{x_j\})=\NESW(\{O_i\},\{x_j\})$. Thus
        \begin{equation}
          \label{eq:GradingD}
          \NESW(\{x_i\}-\{O_i\},(\xUR - \Os)\cap D)=0.
        \end{equation}
        
        Next consider $O_j\in C$, so that its corresponding $X_k$
        which lies in the same column also is in $C$. 
        If the vertical arc connecting $O_j$ to $X_k$ does not cross
        $K_i$, we have that
        \begin{equation*}
          \NESW(\{x_i\},\{O_j\})=\NESW(\{x_i\},\{x_k\})
          \quad\text{and}\quad
          \NESW(\{O_i\},\{O_j\})=\NESW(\{O_i\},\{x_k\}).
        \end{equation*}
        Otherwise both equalities are off by $\pm 1/2$. A more careful
        look at the orientation of the
        horizontal and vertical arcs and our conventions on the
        crossing type reveals that in fact
        \begin{equation}
        \label{eq:GradingC}
        \NESW(\{x_i\}-\{O_i\},(\xUR-\Os)\cap C)
        =
        \#\{\text{negative crossings on $K_i$}\}
        -\#\{\text{positive crossings on $K_i$}\}.
        \end{equation}

        Finally, we claim that
        \begin{equation}
        \label{eq:GradingAB}
        \NESW(\{x_i\}-\{O_i\},(\xUR-\Os)\cap (A\cup B)) =
        -\#\{\text{downward cusps among $\{X_i, O_i\}$}\}.
        \end{equation}
        This follows from an analysis of the eight cases:
        whether the $O_j\in A$ is above or below $X_i$, whether
        the $X_k\in B$ is above or below $O_i$, and whether
        $O_i$ is to the left or to the right of $X_i$. These eight
        cases are illustrated in Figure~\ref{fig:CuspCount}.
\begin{figure}
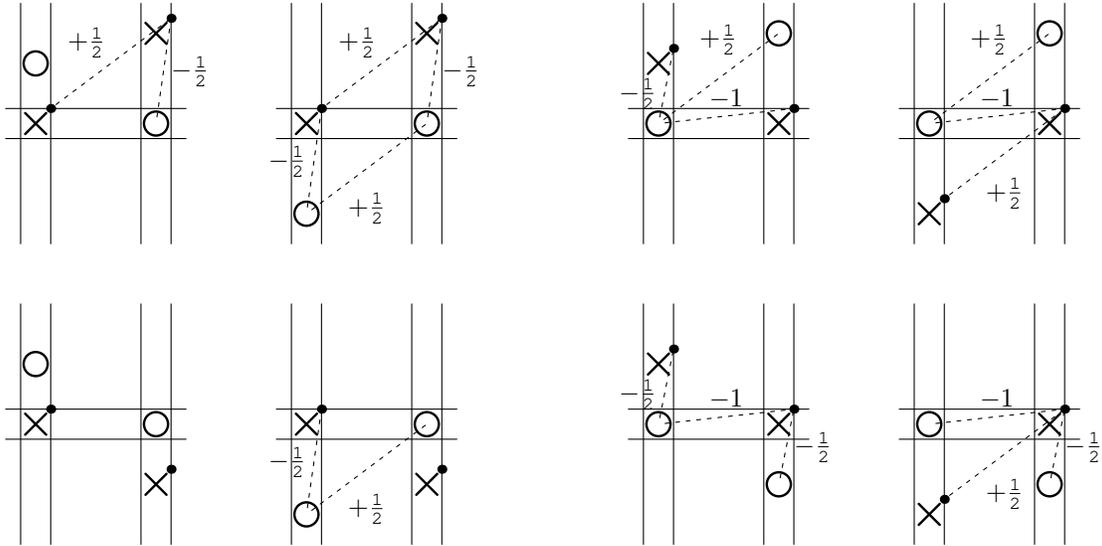

\begin{center}
\graph{cusp-count.10}
\end{center}
\caption {By considering the eight illustrated possibilities, we
  verify Equation~\eqref{eq:GradingAB}.  The relevant contributions to
  $C_i$ are indicated by dashed lines, labelled by the value of the
  contribution.  In each case the sum of the contributions equals the
  number of downward-oriented cusps in the front projection.}
\label{fig:CuspCount}
\end{figure}

Equation~\eqref{eq:GradingPerArc} now follows by adding up
Equations~\eqref{eq:GradingAB}, \eqref{eq:GradingC}, and
\eqref{eq:GradingD}. Equation~\eqref{eq:MasUR} follows from
Equation~\eqref{eq:GradingPerArc} by adding up the contributions of
each horizontal arc, and then adding one (as in
Equation~\eqref{eq:DefMaslov}). Equation~\eqref{eq:MasLL} follows from
a similar
analysis, except that in this case, the contribution from the regions
in $A$ and $B$ are different, and so we replace
Equation~\eqref{eq:GradingAB} by the following:
        \begin{equation}
        \NESW(\{x_i\}-\{O_i\},(\xLL-\Os)\cap (A\cup B)) =
        -\#\{\text{upward cusps among $\{X_i, O_i\}$}\}.\qedhere
        \end{equation}
\end{proof}

\begin{lemma}
        \label{lemma:LegendrianStabilization} Let $G$ be a grid
        diagram and $H$ a stabilization of $G$ of type $\XNW$ or
        $\XSE$.  Then the destabilization map from $\Cm(H)$ to the
        mapping cone of $$U_1-U_2\colon \Cm(G)[U_1]\longrightarrow
        \Cm(G)[U_1]$$
        carries the two elements $\zUR(H)$ and $\zLL(H)$
        to $\pm\zUR(G)$ and $\pm\zLL(G)$, thought of as an element of
        $\Cm(G)\subset \Right\subset C'$, in the notation of
        Section~\ref{subsec:Stabilization} (i.e., where the three
        markings involved are $X_1$, $O_1$, and~$X_2$).
\end{lemma}

\begin{proof}
        This follows from a case analysis of the
        stabilizations.

\begin{figure}
\begin{center}
\graph{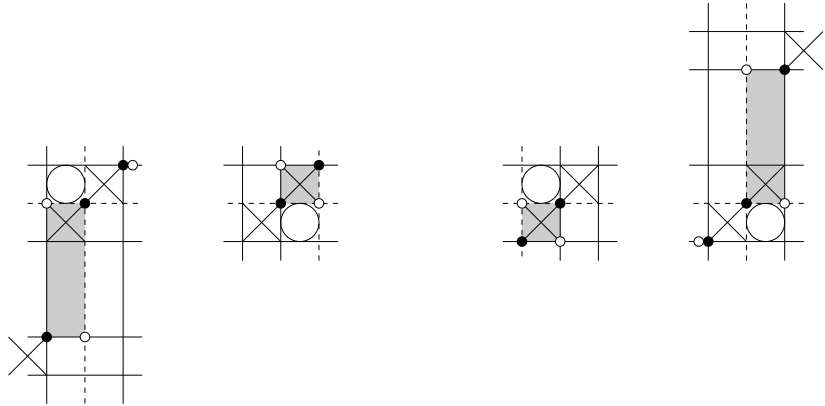}
\end{center}
\caption {{\bf Invariance of the Legendrian element.}
        We have illustrated a local picture of a stabilized
        diagram. The black dots represent the canonical elements for
        the stabilized diagram, the white ones represent the elements
        afterwards, i.e., they are destabilized at the pair of dotted
        lines. The left two diagrams show~$\xUR'$, while the right two
        show~$\xLL'$.  In all four cases the destabilization map counts
        the shaded rectangle, a region with complexity~$2$ of type~$R$.}
\label{fig:LegInvDestab}
\end{figure}

        Let $\xUR{}'$ be the canonical generator for the stabilized
        diagram, and $\xUR$ the corresponding generator in the
        destabilized diagram, considered as an element of~$\Right$.
        
        In each case, we claim that there is only one non-zero term in
        $F(\xUR{}')$, since there is exactly one domain in
        $\pi^F(\xUR{}',\y)$ which does not contain any of the $X_i$
        with $i\neq 1$.  As in Lemma~\ref{lemma:IsACycle}, this
        follows by considering the upper right corner of any such
        domain, which must be~$X_1$. Moreover, we claim that for this
        domain, the image point~$\y$ corresponds to the canonical
        generator $\xUR{}$ for the destabilized picture. We have illustrated
        the two cases on the left of Figure~\ref{fig:LegInvDestab}.
        Because there
        is only one domain, the destabilization map carries $\zUR(H)$
        to $\pm\zUR(G)$ as desired.  As sketched on the right of
        Figure~\ref{fig:LegInvDestab}, a similar argument works for
        $\xLL$.
\end{proof}

\begin{lemma}
  \label{lemma:Commutation}
        Let $G$ and $H$ be two grid diagrams which differ by a
        commutation move. Then under the map
        $F_{\beta\gamma}\colon \Cm(G) \longrightarrow \Cm(H)$, the
        image of the cycle $\zUR(G)$ is $\pm\zUR(H)$ and the image of
        $\zLL(G)$ is $\pm\zLL(H)$.
\end{lemma}

\begin{proof}
        We argue that there is exactly one $\y$ and one pentagon in
        $\Pent(\xUR(G),\y)$ which does not contain some $X_i$ in its
        interior, and that is the one which connects $\xUR(G)$ to
        $\xUR(H)$. This can be seen from an argument like that in
        Lemma~\ref{lemma:LegendrianStabilization}: there is a pentagon
        taking $\xUR(G)$ to $\xUR(H)$, as illustrated in
        Figure~\ref{fig:FindPentagons}.  Suppose
        that $\psi$ is any other pentagon, and consider its upper
        right corner, which is at some
        $c\in\xUR(G)$. The subsquare just to the lower left of this
        $c$ contains some $X_i$ with $i\neq 1$ or $2$.  The argument
        for $\xLL$ is similar.
\begin{figure}
\begin{center}
\begin{picture}(0,0)%
\graph{CommuteNat.pstex}%
\end{picture}%
\setlength{\unitlength}{1579sp}%
\begingroup\makeatletter\ifx\SetFigFont\undefined%
\gdef\SetFigFont#1#2#3#4#5{%
  \reset@font\fontsize{#1}{#2pt}%
  \fontfamily{#3}\fontseries{#4}\fontshape{#5}%
  \selectfont}%
\fi\endgroup%
\begin{picture}(2424,8886)(1189,-7294)
\put(2851,1289){\makebox(0,0)[lb]{\smash{{\SetFigFont{9}{10.8}{\rmdefault}{\mddefault}{\updefault}{\color[rgb]{0,0,0}$\beta$}%
}}}}
\put(1651,1289){\makebox(0,0)[lb]{\smash{{\SetFigFont{9}{10.8}{\rmdefault}{\mddefault}{\updefault}{\color[rgb]{0,0,0}$\gamma$}%
}}}}
\put(2701, 89){\makebox(0,0)[lb]{\smash{{\SetFigFont{9}{10.8}{\rmdefault}{\mddefault}{\updefault}{\color[rgb]{0,0,0}$b$}%
}}}}
\put(1651,-4111){\makebox(0,0)[lb]{\smash{{\SetFigFont{9}{10.8}{\rmdefault}{\mddefault}{\updefault}{\color[rgb]{0,0,0}$a$}%
}}}}
\end{picture}%
\end{center}
\caption {{\bf Small pentagon.}
\label{fig:FindPentagons}
        The dark circles represent the canonical generator $\xUR$
        for $G$, the diagram involving $\beta$, while the empty
        circles represent the canonical generator $\xUR{}'$ for
        $H$, the diagram involving $\gamma$. The shaded pentagon
        represents the map carrying $\xUR$ to $\xUR{}'$.}
\end{figure}
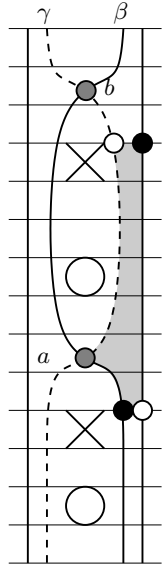
\end{proof}

\begin{proof}[Proof of Theorem~\ref{thm:LegendrianInvariant}, up to signs.]
The calculations of the Maslov and Alexander gradings of~$\zUR$
are Lemmas~\ref{lemma:Maslov} and~\ref{lemma:AlexanderGrading} respectively.
The Alexander grading is given by 
\begin{align*}
A(\zUR) &= -\writhe(K)-\#\{\text{downward-oriented cusps}\} \\
&= \writhe(\Mirror(K))-\#\{\text{downward-oriented cusps}\}, \\
A(\zLL) &=\writhe(\Mirror(K))-\#\{\text{upward-oriented cusps}\}.
\end{align*}
Comparing these with the standard descriptions of $\TB$ and
$\Rotation$ from the knot projection (Equations~\eqref{eq:TBformula}
and~\eqref{eq:Rformula} respectively), we obtained the stated
formulas for the bigradings of $\zUR$ and $\zLL$.

In view of Proposition~\ref{prop:LegendrianGrid}, invariance under
Legendrian isotopies (up to sign) follows from
Lemmas~\ref{lemma:LegendrianStabilization}
and~\ref{lemma:Commutation}.
\end{proof}

\subsection{Properties of the Legendrian invariant}

We now turn to the properties of $\UR$ and $\LL$ stated in the
introduction.

\begin{proof}[Proof of Proposition~\ref{prop:Symmetries}, up to signs.]
The fact that the two invariants are permuted under orientation
reversal follows from the symmetry of the torus which is given by
reflection through the $x=-y$ axis. More precisely, it is easy to
see that if $G$ is a grid diagram representing a Legendrian
knot~$\oLegK$ and $H$ is the grid diagram obtained
from~$G$ by reflecting through this axis, then the reflection map from
$\S(G)$ to $\S(H)$, which takes $\xUR(G)$ to $\xLL(H)$ and $\xLL(G)$
to $\xUR(H)$, induces
an isomorphism of complexes
$$\Phi\colon \Cm(G) \longrightarrow \Cm(H),$$
with $\Phi(\zUR(G))=\pm\zLL(H)$ and $\Phi(\zLL(G))=\pm\zUR(H)$. Moreover,
by Lemma~\ref{lem:GridSymmetries}, the Legendrian knot specified
by~$H$ is~$-\oLegK$.

Similarly, rotation of a grid diagram~$G$ by~$180\deg$ induces an
isomorphism of chain complexes which permutes the two canonical cycles
while taking~$\oLegK$ to its Legendrian mirror~$\mu(\oLegK)$.
\end{proof}

\begin{figure}
\begin{center}
\graph{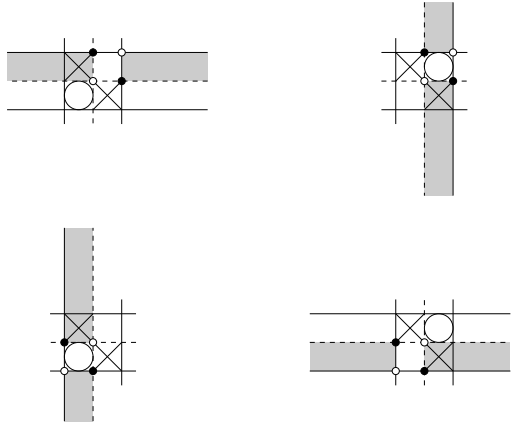}
\end{center}
\caption {{\bf Behaviour under positive and negative destabilization.}
\label{fig:PosNegDestabilization}
        The black generators represent the canonical elements for the
        stabilized diagram, while the white ones represent canonical
        elements in the destabilized one. The left two diagrams are
        positive destabilization (of type $\XNE$) and the right two
        are negative destabilization (of type $\XSW$).  On the top row,
        we consider~$\xUR$, while on the bottom row, we consider~$\xLL$.}
\end{figure}

\begin{proof}[Proof of Theorem~\ref{thm:StabilizationTheorem}, up to signs.]
  Behaviour under positive and negative destabilization is illustrated
  in Figure~\ref{fig:PosNegDestabilization}.  The illustrated domains,
  all of complexity~$2$ and type~$R$, are as before the only domains
  of type~$F$ starting at~$\xUR$ or~$\xLL$ which do not contain
  any other~$X_i$ in their interior.  In all cases the rectangle
  connects $\x^{\pm}$ to the
  canonical generator for the destabilized picture.  For the positive
  destabilization of~$\xUR$ and the negative destabilization
  of~$\xLL$, the domain contains the other~$O$ in the row
  containing~$X_1$ (and so the map on homology is multiplication
  by~$U$), while for the other two cases the domain contains
  only~$X_1$.  Therefore the induced map~$\Phi$ on the chain complex
  behaves as stated, up to a sign.
\end{proof}

In order to prove the invariant is non-zero
(Theorem~\ref{thm:NonVanishing}), we will look at
yet another complex.

\begin{definition}
  Let the complex~$C^!(G)$ be the tensor product of~$\CKm(G)$ with
  Laurent polynomials in the $\{U_i\}_{i=1}^n$, i.e.,
  $\Z[U_1,U_1^{-1},\dots,U_n,U_n^{-1}]$.  Let $H^!(G)$ be the homology
  of~$C^!(G)$, which is a module over $\Z[U,U^{-1}]$.
\end{definition}

We will find it useful to use the following:

\begin{definition}
  A {\em domain} $p$ from $\x$ to $\y$ is a two-chain in $\Torus$
  whose boundary $\partial p$ is a path from $\x$ to $\y$. We denote
  the set of domains from $\x$ to $\y$ by $\pi(\x,\y)$.  Given
  $\phi\in\pi(\x,\y)$, let $X_j(\phi)$ resp. $O_j(\phi)$ denote the
  local multiplicity of $\phi$ at $X_j$ resp. $O_j$.
\end{definition}

\begin{lemma}
  \label{lem:HomologyShriek}
  $H^!(G)$ is isomorphic to $\Z[U,U^{-1}]$, and it is generated
  by $[\xUR(G)]$ or
  $[\xLL(G)]$.
\end{lemma}

\begin{proof}
  The chain complexes~$C^!$ for different positions of the $O$'s are
  isomorphic, as follows.  Suppose that $G$ is a grid diagram with
  grid number $m$, with two alternate possible placement of the $O$'s:
  $\Os$ and $\Os'$, and let $C^!(G)$ and $C^!(G')$ be the two possible
  corresponding chain complexes. We will construct an isomorphism
  $\Phi\colon C^!(G) \longrightarrow C^!(G')$.
  
  Fix $\x_0\in\S(G)$ arbitrarily, let $\z_0$ and $\z_0'$ be the
  corresponding
  elements of $C^!(G)$ and $C^!(G')$, and declare $\Phi(\z_0)=\z_0'$.
  Given any generator of $C^!(G)$ {\em as a $\Z[U_1,U_1^{-1}]$-module}, which
  has the form of $U_2^{n_2}\cdots U_m^n\cm \y$, it is easy
  to see that there is a unique $\phi\in\pi(\x,\y)$ with $X_j(\phi)=0$
  for all $j$, and $O_j(\phi)=n_j$ for $j=2,\dots,m$.  (We use here
  the fact that $G$ represents
  a knot, rather than a link with more than one component.)
  
  We then define $\Phi(U_2 \cdots U_m\cm \y)=U_2^{O_2'(\phi)}
  \cdots U_m^{O_m'(\phi)}\cm \y$. This is easily seen to be an
  isomorphism of chain complexes of $\Z[U,U^{-1}]$-modules, where
  we take $U=U_1$.
  
  In particular the complex~$C^!(G)$ can be identified with the
  complex $C^!(G')$, where $G'$ is the grid diagram with one~$O$
  directly beneath each~$X$; note that $G'$ is a diagram for the
  unknot.  Furthermore, note that $[\zUR(G')]$ in $H^!(G')$ is
  invariant up to units under \emph{all} destabilizations, by
  Theorems~\ref{thm:StabilizationTheorem} (since multiplication
  by~$U_i$ is now invertible).  Therefore we can simplify~$G'$ until
  we get a $2\times 2$ grid diagram~$G_0$ representing a trivial
  unknot, where an elementary calculation shows $H^!(G_0)$ is rank
  one, generated by~$[\zUR(G_0)]$.  Again, similar computations work
  for~$\xLL(G)$.
\end{proof}

\begin{proof}[Proof of Theorem~\ref{thm:NonVanishing}]
  By Lemma~\ref{lem:HomologyShriek}, the map from $\HFKm(G)$ to $H^!(G)$
  takes $\UR(G)$ to a generator of $H^!(G)$.  In particular
  $U^m\UR(G) \ne 0$ for any non-negative integer~$m$.
\end{proof}

We can now complete Definition~\ref{def:Canonical}.

\begin{definition}
  \label{def:CanonicalWithSigns}
  The signs in the definition of~$\zLL(G)$ and~$\zUR(G)$ in
  Definition~\ref{def:Canonical} are chosen so that
  $\zUR(G)=\xUR(G)$, and there is an~$m\in\Z$ so that
  $[\zLL(G)]=U^m[\zUR(G)]$ when thought of as
  elements of $H^!(G)$.
\end{definition}

\begin{proof}[Proof of Theorem~\ref{thm:LegendrianInvariant}, with signs]
With the above definition, we can see that in
Lemmas~\ref{lemma:LegendrianStabilization}
and~\ref{lemma:Commutation}, either $\zUR(G)$ and $\zLL(G)$ map to
$\zUR(H)$ and $\zLL(H)$ respectively, or they map to $-\zUR(H)$ and
$-\zLL(H)$, depending on whether the chosen generator of $H^!(G)$ maps
to the chosen generator of $H^!(H)$ or its negative.  By negating
the chain map~$\Phi$ from $\Cm(G)$ to $\Cm(H)$ if
necessary, we can make it map $\zUR(G)$ to $\zUR(H)$ and $\zLL(G)$ to
$\zLL(H)$, with correct signs.
\end{proof}

 Similarly we can fix the signs in
Proposition~\ref{prop:Symmetries} and
Theorem~\ref{thm:StabilizationTheorem}.


\section{The case of links}
\label{sec:Links}

Most of the discussion from the earlier parts of this paper have
rather straightforward generalizations to the case of links, as we
will now show.

Let $\oLegL$ be an oriented Legendrian link with $\ell$
components. The contact distribution
$\xi$ determines a complex line bundle over $S^3$ equipped with a
trivialization in a neighborhood of $\oLegL$. The Euler
class of this line bundle relative to its trivialization on the boundary
gives an element of $H^2(S^3,\oLegL;\Z)\cong H_1(S^3-L;\Z)\cong
\Z^\ell$, which plays the role of the rotation number. 
More concretely, we obtain $\ell$ integers, $\{\Rotation_i\}_{i=1}^\ell$,
determined by the property that if $F$ is any surface in $S^3$ whose
boundary lies on~$\oLegL$, then
$$\langle e(\xi|_{\partial F},\oLegL'), [F]\rangle
= \sum_{i=1}^\ell \Rotation_i\cm \#(F\cap m_i),$$
where $\{m_i\}_{i=1}^\ell$ are the meridians for the components of $L$.
Similarly, the Thurston\hyp Bennequin framing gives an element of $H_1(S^3-L;\Z)$,
which we can write as 
$$\sum_{i=1}^\ell \TB_i\cm m_i.$$
That is, $\TB_i$ is the linking number of
the Legendrian push-off of $\oLegL$ with the $i\th$ component
of~$\oLegL$.  Note that $\TB_i$ depends on the orientation
of~$\oLegL$.
One alternate (equivalent) data to $\TB_i$ is the Thurston-Bennequin
invariant of
the $i\th$ component considered as a knot by itself; that choice is
less convenient for us.

In terms of the front projection $\Pi=\bigcup_{i=1}^\ell\Pi_i$, we have that
\begin{align}
\label{eq:TBProj}
\TB_i(\oLegL)&=\writhe(\Pi_i)+\lk(\Pi_i,\Pi-\Pi_i)
-\OneHalf\#\{\text{cusps in $\Pi_i$}\}\\
\label{eq:RProj}
\Rotation_i(\oLegL)&=\OneHalf \Big(\#\{\text{downward-oriented cusps in $\Pi_i$}\} -
\#\{\text{upward-oriented cusps in $\Pi_i$}\}\Big).
\end{align}

Let $G$ be a grid diagram for $\oLegL$. We can, of course,
define the chain complex $\Cm(G)$ as before, but in fact, this complex
also has a refinement. Specifically, we can consider the
Alexander multi-grading, which is a function $A \colon
\S\longrightarrow \Z^\ell$ defined as follows. We partition $\Os$ as
$\bigcup_{i=1}^\ell \Os_i$ and $\Xs$ as $\bigcup_{i=1}^\ell\Xs_i$,
where $\Os_i$
(respectively $\Xs_i$) denotes the set of~$O_j$ (respectively $X_j$)
corresponding to
the $i\th$ component of the link.  Let $n_i= \#\Os_i$. Now we can
define
$A(\x)=(A_1(\x),\ldots,A_\ell(\x))$, where
\begin{equation}
  \label{eq:AlexanderFormulaTwo}
A_i(\x) = \NESW\Big(\x-\OneHalf(\Xs+\Os),\Xs_i-\Os_i\Big)- \Bigl(\frac{n_i-1}{2}\Bigr).
\end{equation}

We can form the chain complex $\CLm(\oLink)$, defined as in
Equation~\eqref{eq:DefCLm}.  We number the variables so that the
first $\ell$ of the $O$'s, $O_1,\dots,O_\ell$, belong to the $\ell$ different
components of the link. Taking the homology of this module, we obtain
a graded module over $\Ring=\Z[U_1,\dots,U_\ell]$
$$\HFLm(L)=\bigoplus_{\substack{d\in\Z\\
    \s\in\Z^\ell}}\HFLm_d(L,\s),$$
where $U_i$ acts as an endomorphism
which is homogeneous of degree~$-2$ for~$d$, degree~$-1$ for the
$i\th$ component of~$\s$, and degree~0 otherwise. This $\Ring$-module
is an oriented link invariant~\cite{MOST} which, when
specialized to coefficients in $\Field$, agrees with link Floer
homology~\cite{Links}. 

Definition~\ref{def:Canonical} readily generalizes to this context,
giving a pair of elements $\zLL, \zUR\in\Cm(G)$, each of which is a
cycle in the associated graded object~$\CLm(G)$.  We have the
following analogue of Theorem~\ref{thm:LegendrianInvariant}.

\begin{theorem} 
        \label{thm:LegendrianInvariantLink}
        For a grid diagram~$G$, let $\oLegL=\oLegL(G)$ be the
        corresponding oriented Legendrian link. Then there are two
        associated cycles $\zUR=\zUR(G)$ and $\zLL=\zLL(G)$, supported
        in gradings
\begin{align}
A_j(\zUR)&=\frac{\TB_j(\oLegK)-\Rotation_j(\oLegK)+1}{2}
  &A_j(\zLL)&=\frac{\TB_j(\oLegK)+\Rotation_j(\oLegK)+1}{2}
 \label{eq:AlexLinkInv}\\
M(\zUR)&=2 \sum_{j=1}^\ell A_j(\zUR)+1-\ell 
  &M(\zLL)&= 2 \sum_{j=1}^\ell A_j(\zLL)+1-\ell,
\label{eq:MaslovLinkInv}
\end{align}
where here $\ell$ denotes the number of components of $\oLegL$.
Moreover, if $G$ and $G'$ are two different grid diagrams which
represent Legendrian isotopic oriented links, then there is 
a quasi-isomorphism of chain complexes
$$\Phi\colon \Cm(G) \longrightarrow \Cm(G')$$
with
\begin{align*}
\Phi(\zUR(G))&=\zUR(G')& \Phi(\zLL(G))&=\zLL(G').
\end{align*}
\end{theorem}
\begin{proof}
Most of this proof is a straightforward generalization of the proof of
Theorem~\ref{thm:LegendrianInvariant}. For example,
Equation~\eqref{eq:MaslovLinkInv} is  a straightforward adaptation of
the argument from Lemma~\ref{lemma:AlexanderGrading}.

Equation~\eqref{eq:AlexLinkInv} follows along the lines of
the proof of Lemma~\ref{lemma:Maslov}, but with a little extra care.
We first argue that
\begin{equation}
  \label{eq:AlexanderForLinks}
 2A_j(\zUR)=-\writhe(K)-\#\{\text{downward cusps on component
 $j$}\}+1. 
\end{equation}
As in Lemma~\ref{lemma:Maslov}, we consider the horizontal arc $K_i$ 
connecting
$X_i$ and $O_i$, and the corresponding partition of the plane into regions
$A$, $B$, $C$, and $D$.  We also let $K_i$ denote the corresponding
horizontal strip in the plane. We then analyze the contributions
to $2\NESW\left(K_i\cap (\xUR-\OneHalf(\Xs+\Os)), (\Xs_j-\Os_j)\right)$
coming from the portion of $(\Xs_j-\Os_j)$ in these four possible regions.

As in Equation~\eqref{eq:GradingD}, we have
\begin{equation}
2\NESW\left(K_i\cap (\xUR-\OneHalf(\Xs+\Os)),D\cap
  (\Xs_j-\Os_j)\right)=0.
\end{equation}
As in Equation~\eqref{eq:GradingC}, provided that $K_i\subset L_j$, we have
\begin{multline}
      2\NESW\left(K_i\cap (\xUR-\OneHalf(\Xs+\Os)),C\cap (\Xs_j-\Os_j)\right)
        = \#\{\text{negative crossings of $L_j$ with $K_i$}\} \\
        -\#\{\text{positive crossings of $L_j$ with $K_i$}\}.
\end{multline}
Finally, we consider the analogue of 
Equation~\eqref{eq:GradingAB}. 
If $K_i$ is an arc on
component~$j$, by considering the eight cases from
Figure~\ref{fig:CuspCount}, we have that
\begin{multline}
        2\NESW\left(K_i\cap (\xUR-\OneHalf(\Xs+\Os)), (A\cup B)\cap
        (\Xs_j-\Os_j)\right)=\\
\begin{aligned}
&1-\#\{\text{downward cusps among $\{X_i, O_i\}$}\} \\
&- (\text{$\OneHalf$ if $\Os\cap K_i$ is above the $X$ in its column}) \\
&+ (\text{$\OneHalf$ if $\Xs\cap K_i$ is below the $O$ in its column}).
\end{aligned}
\end{multline}
Otherwise, this contribution is zero.  We deduce
Equation~\eqref{eq:AlexanderForLinks} by adding up half the local
contributions calculated above and subtracting
$\frac{n_i-1}{2}$, as in Equation~\eqref{eq:AlexanderFormulaTwo}. The statement for $A_j(\xUR)$ follows now from that equation,
together with Equations~\eqref{eq:TBProj} and~\eqref{eq:RProj}.

The definition of $\xUR$ and $\xLL$ is as before.  The proof of
invariance under the chain map up to signs follows identically.
The signs in the definition of $\zUR$ and $\zLL$ are fixed later, in
Definition~\ref{def:CanonicalWithSignsLinks}.
\end{proof}

Let $\UR(\oLegL)$ and $\LL(\oLegL)$ denote the homology classes in
$\HFLm(L)$ of $\zUR(G)$ and $\zLL(G)$.  These are
Legendrian invariants of $\oLegL$.  Behaviour under orientation
reversal of all components simultaneously and Legendrian mirror is the
same as in Proposition~\ref{prop:Symmetries}. The analogue of
Theorem~\ref{thm:StabilizationTheorem} is the following.

\begin{theorem}
        \label{thm:StabilizationTheoremLink} Let $\oLegL$ be an oriented
        Legendrian link, and $\oLegL^-$ (respectively $\oLegL^+$) be 
        oriented Legendrian links obtained as a single negative
        (respectively positive) stabilization of $\oLegL$ on the $j\th$
        component. Then there is a
        quasi-isomorphism 
        \begin{align*} 
                \Phi^-\colon C(\oLegL) &\longrightarrow C(\oLegL^-) \\
                \Phi^+\colon C(\oLegL) &\longrightarrow C(\oLegL^+) 
        \end{align*} 
        under which
        \begin{align*}
                 \Phi^-(\UR(\oLegL))&=\UR(\oLegL^-) &
                 U_j\cm \Phi^-(\LL(\oLegL))&= \LL(\oLegL^-) \\
                 U_j\cm \Phi^+(\UR(\oLegL))&=\UR(\oLegL^+) &
                 \Phi^+(\LL(\oLegL))&= \LL(\oLegL^+).
        \end{align*}
\end{theorem}

\begin{proof}
  This is a straightforward generalization of the proof of
  Theorem~\ref{thm:StabilizationTheorem}.
\end{proof}

Similarly, we have the following generalization of Theorem~\ref{thm:NonVanishing}:

\begin{theorem}
  \label{thm:LegendrianLink}
  For a Legendrian link $\oLegL$, the homology classes $\UR(\oLegL)$
  and $\LL(\oLegL)$ are non-trivial; and indeed, they are not
  $U_i$-torsion for any of the $U_i$.
\end{theorem}

\begin{proof}
  A little more care is needed than in the case of knots: $H^!(G)$ now
  does depend on the
  placement of the $O_i$. However, if we consider the new complex
  $C'(G)=\Cm(G)/\{U_i=1\}_{i=1}^m$, then the homology of this complex
  $H'(G)$ obviously no longer depends on the placement of the $O_i$
  (since they play no role in the $C'(G)$).  Thus we can move
  placements of the $\Os$ to realize the unknot.  But we already know
  by Theorem~\ref{thm:NonVanishing} that the invariants for the unknot
  do not vanish in $\HFKm$, so they do not vanish in $H'(G)$,
  so both $\UR(\oLegL)$ and $\LL(\oLegL)$ represent non-trivial homology
  classes in $H'(G)$. It follows readily that the two classes
  are not $U_i$-torsion for any of the $U_i$.
\end{proof}

As a final point in this section, we turn to the signs entering the
definition of $\UR(L)$ and $\LL(L)$, generalizing
Definition~\ref{def:CanonicalWithSigns} to the case of links.

Let $C''(G)=C(G)/\{U_1=\dots=U_\ell=1\}$, where as
usual we number the $U_i$ so that the first $\ell$ correspond to the
$\ell$ distinct components of the link. Although the complex no longer
inherits neither an Alexander of Maslov gradings, $C''$ retains an
integral grading given by $N=M-2A+\ell-1$.

\begin{lemma}
  \label{lem:RelSigns}
  We have that $H_*(C''(G))\cong H_*(T^{\ell-1})$,
  and the elements $[\xUR(G)]$ and $[\xLL(G)]$ generate
  the zero-dimensional part.
\end{lemma}

\begin{proof}
  Let $C'(G)$ denote the complex studied in the proof of
  Theorem~\ref{thm:LegendrianLink}.  As in that proof, its homology is
  independent of the placement of the $O_i$, so it is isomorphic to
  the value for the unknot, from which it is easy to see that
  $H_*(C'(G))\cong H_*(T^{m-1})$ with an appropriate overall grading
  shift.  (Remember that $H_*(T^{m-1})$ is a $(m-1)$-fold tensor
  product of a space $V'$ with itself, where $V' \cong \Z \oplus \Z$
  with one generator in degree~0 and one in degree~1.)
  
  Note that $C'(G)=C''(G)/\{U_{\ell+1}=\dots=U_m=1\}$. Moreover, in
  $C'(G)$, each $U_i$ is homotopy equivalent to multiplication by~$1$.
  (This follows from the fact that in $\CLm$, multiplication by $U_i$
  and multiplication by $U_j$ are chain homotopic if $O_i$ and $O_j$
  correspond to the same component of the link,
  see \cite[Lemma~2.9]{MOST}, together with the fact that we have set
  $U_i=1$ for $i=1,\dots,\ell$ in $C''(G)$.)  Thus, it follows
  that $H_*(C'(G))\cong H_*(C''(G))\otimes H_*(T^{m-\ell})$. Combining
  this with the fact that $H_*(C'(G))\cong H_*(T^{m-1})$, we can
  conclude that $H_*(C''(G))\cong H_*(T^{\ell-1})$.

  From the proof of Theorem~\ref{thm:LegendrianLink},
  $[\xUR(G)]$ and $[\xLL(G)]$ are primitive elements in $H_*(C''(G))$.
  From Equation~\eqref{eq:MaslovLinkInv}, we conclude that they both
  are supported in degree zero.
\end{proof}

In view of Lemma~\ref{lem:RelSigns}, can generalize
Definition~\ref{def:CanonicalWithSigns} to the case of links, as
follows.

\begin{definition}
  \label{def:CanonicalWithSignsLinks}
  The signs in the definition of $\zUR$ and $\zLL$ are chosen
  so that $[\zUR(G)]$ and $[\zLL(G)]$ represent the same class in the
  homology of
  $C''(G)=C(G)/\{U_1=\dots=U_\ell=1\}$.
\end{definition}

This is clearly compatible with
Definition~\ref{def:CanonicalWithSigns} in the case
where $\ell=1$, and as before we can now fix the signs in
Theorems~\ref{thm:LegendrianInvariantLink}
and~\ref{thm:StabilizationTheoremLink}.


\section{Examples: The knot $5_2$ and the link $6^2_3$}
\label{sec:examples}

To demonstrate that our invariant is non-trivial, we will use it to
distinguish different Legendrian representatives of the knot~$5_2$ and
the link~$6^2_3$.  Throughout this section we will work over $\Field$.
\begin{example}
  \label{ex:Chekanov}
  The unoriented Legendrian knots~$\LegK_1$ and~$\LegK_2$ with front
  projections
  \[
  \LegK_1 = \mfigb{grids.5}\quad\text{and}\quad\LegK_2 = \mfigb{grids.15},
  \]
  both of topological type $5_2$ and having $\TB=1$ and
  $\Rotation=0$, are not Legendrian isotopic.
\end{example}

This example was first found by Chekanov~\cite{Chekanov}.

\begin{proof}
  Reflecting through the vertical axis takes both front projections
  shown above to themselves with orientation reversed.  As in
  Lemma~\ref{lem:GridSymmetries}, reflecting the front through the
  vertical axis corresponds to rotating the Legendrian knot by
  $180\deg$ around the $z$-axis, which is an isotopy.
  Therefore the two oriented Legendrian knots are both
  isotopic to their reverse and it suffices to prove the result for
  one orientation.

  Grid diagrams $G_1$, $G_2$ for the two Legendrian knots are shown in
  Figure~\ref{fig:fiveTwo}.
  On~$G_1$, we have indicated the generator~$\zUR$.  From the diagram
  it is straightforward to check that $\zUR$ is isolated in the chain
  complex: for all generators~$\x$, there are no rectangles in
  $\EmptyRect(\x,\zUR)$ with empty intersection with $\Xs$ and $\Os$.
  By the symmetry of the diagram,
  the same is true for~$\zLL$.  It follows that $\UR(\oLegK_1) \ne
  \LL(\oLegK_1)$.
\begin{figure}
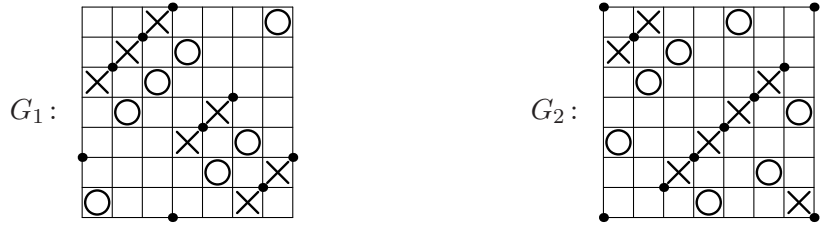

  $G_1\colon\mfigb{grids.1}\hspace{3cm}G_2\colon\mfigb{grids.13}$
  \caption{Grid diagrams for two different Legendrian representatives
    of the knot~$5_2$.  The left diagram shows the generator $\zUR$.
    The right diagram shows a generator~$\y$ so that $\partial\y =
    \zUR + \zLL$.}
  \label{fig:fiveTwo}
\end{figure}

On the other hand, consider the generator~$\y$ shown on~$G_2$ in
Figure~\ref{fig:fiveTwo}.  An elementary check shows
that $\partial(\y) = \zUR + \zLL$; therefore $\UR(\oLegK_2) =
\LL(\oLegK_2)$.  But if $\oLegK_1$ and $\oLegK_2$ were Legendrian
isotopic (with either orientation), by
Theorem~\ref{thm:LegendrianInvariant} there would be an isomorphism from
$\HFKm(\oLegK_1)$ to $\HFKm(\oLegK_2)$ taking~$\UR(\oLegK_1)$ to~$\UR(\oLegK_2)$ 
and~$\LL(\oLegK_1)$
to~$\LL(\oLegK_2)$, a contradiction.
\end{proof}

\begin{remark}
  \label{rem:InvarianceLimits}
  This example shows that the transverse invariant
  $\TransElt(\TransK)$ is only an invariant up to quasi\hyp
  isomorphisms (and not as an element of an abstract homology group).
  In particular, a sequence of elementary moves from a grid
  diagram~$G$ back to itself induces a quasi\hyp isomorphism that need
  not be the identity.  Indeed, an elementary calculation shows that
  the knots~$\oLegK_1$ and~$\oLegK_2$ become Legendrian isotopic after
  either one positive or one negative stabilization.  If we follow a
  path of grid diagrams that takes~$G_1$ to~$G_2$ via a positive
  stabilization and destabilization, followed by a path that
  takes~$G_2$ to~$G_1$ via a negative stabilization and
  destabilization, we first take $\LL(\oLegK_1)$ to
  $\LL(\oLegK_2)$, which is equal to $\UR(\oLegK_2)$, which we take to
  $\UR(\oLegK_2)$.  Since $\LL(\oLegK_1) \ne\UR(\oLegK_2)$, this
  sequence of grid moves induces a non\hyp trivial automorphism
  of~$\HFKm(\oLegK_1)$.
\end{remark}

\begin{example}
  The oriented Legendrian links $\oLegL_3$ and~$\oLegL_4$ with front
  projections
  \[
  \oLegL_3 = \mfigb{grids.80}
  \quad\text{and}\quad
  \oLegL_4 = \mfigb{grids.70},
  \]
  both of topological type $6^2_3$ and having $\TB_i=1$ and
  $\Rotation_i=0$ for $i=1,2$,
  are not Legendrian isotopic.
\end{example}

\begin{proof}
  Grid diagrams $G_3$ and $G_4$ for these two links are shown in
  Figure~\ref{fig:twistedWhitehead}.  It is easy to verify that
  $\zUR(G_3)$ and $\zLL(G_3)$ are isolated in the chain complex, as
  for $G_1$.  It follows that
  $\UR(\oLegL_3)\ne \LL(\oLegL_3)$.  On the other hand,
  we have indicated a generator $\y$ on
  $G_4$ so that $\partial(\y) = \zUR + \zLL$, so
  $\UR(\oLegL_4) = \LL(\oLegL_4)$.
\end{proof}

\begin{figure}
  \[
  G_3\colon\mfigb{grids.82}\hspace{2.5cm}G_4\colon\mfigb{grids.74}
  \]
  \caption{Grid diagrams $G_3$ and $G_4$ for $\oLegL_1$ and
    $\oLegL_2$.  The left diagram shows $\zUR(G_3)$.  The right
    diagram shows a generator $\y$ with $\partial \y = \zUR + \zLL$.}
  \label{fig:twistedWhitehead}
\end{figure}


\appendix
\section{On $\tau$}
\label{sec:OnTau}

We found it convenient to work with $\HFKm$ in the present paper.
In particular, in the introduction, we gave a definition of $\tau$ 
which refers to $\HFKm$, i.e., it is defined in terms of the 
associated graded object for the Alexander filtration. 
By contrast, the usual definition of $\tau$ refers to the filtered
complex $\CFa$~\cite{4BallGenus}.
For completeness, we repeat this definition, putting it in terms of
grid diagrams.

Consider the chain complex $\Cm(G)$ as in Equation~\eqref{eq:DefCm}.
This chain complex has a $\Z$ filtration, induced by the Alexander filtration,
whose associated graded object is $\CKm(G)$ considered throughout most of this
paper. 

We will need two constructions, as follows. If we set $U_1=0$, we
obtain a new chain complex $\Ca(G)$, filtered by sub-complexes
${\widehat\Filt}(K,s)$, which are generated by elements with Alexander grading
$\leq s$.  Following~\cite{4BallGenus}, we define
$$\tau(K)=\min\,\{\,s\mid {\widehat \Filt}(K,s)\longrightarrow \HFa(S^3)~\text{is non-trivial}\,\}.$$

More symmetrically, we can consider the filtration
${\widetilde\Filt}(K,s)$ on
$${\widetilde C}(G)={\Cm(G)}/{(U_1=\dots=U_n=0)},$$
where ${\widetilde\Filt}(K,s)$, once again, is generated by those intersection
points $\x$ with $A(\x)\leq s$. This, too, could be used to
calculate~$\tau$: Define
$$\widetilde\tau(K)=\min\,\{\,s\mid H_*({\widetilde \Filt}(K,s))\longrightarrow
  H_*({\widetilde C}(G))~\text{is non-trivial}\,\}.$$

\begin{lemma}
  \label{lemma:TauTilde} $\widetilde\tau(K) = \tau(K)-n+1$.
\end{lemma}

\begin{proof}
  The chain complex ${\widetilde C}(G)$ is filtered quasi-isomorphic
  to the filtered mapping cone of an iterated mapping cylinder
  $$\left(\frac{C(G)}{U_1=0}\right)\otimes_{\Ring}\left(\bigotimes_{i\neq
      1} \Ring[-1,-1] \xrightarrow{U_i}\Ring\right),$$
  where
  $\Ring=\Z[U_1,\dots,U_n]$, and $\Ring[-1,-1]$ denotes $\Ring$ with
  a shift in bigrading so that $1$ has both Alexander filtration and Maslov
  grading of $-1$ (see \cite[Lemma~{2.1}]{MOST}). Since
  multiplication by $U_i$ on $C$ is filtered
  chain homotopic to multiplication by $U_{1}$ (which we have set
  equal to zero in the quotient complex), the above mapping cylinder
  is in fact filtered quasi-isomorphic to
  $$\left(\frac{C}{U_1=0}\right)\otimes V^{n-1}, $$
  where $V$ is a free Abelian group generated by two elements,
  one with Maslov and Alexander bigrading $(0,0)$, and another with 
  Maslov and Alexander bigrading $(-1,-1)$.  In particular, if we
  always choose the second generator of~$V$, we get a sub-complex
  of~$\widetilde C(G)$ which is isomorphic to~$\widehat C(G)$ with
  both gradings shifted by $-n+1$.  The element with minimal Alexander
  filtration mapping non-trivially to $H_*(\widetilde{C}(G))$ will
  always live in this sub-complex.
\end{proof}

To fit the definition of~$\tau$ above into the discussion from the
introduction, define
  $$\tau'(K)=\max\,\{\,s\in \Z \mid\exists\xi\in \HFKm(K,s)~\text{such
    that}~\forall d\geq 0,~U^d\xi\neq0\,\}.$$

\begin{lemma} $\tau'(K) = \tau(\Mirror(K))$.
\end{lemma}

\begin{proof}
  Let $C'(G)=\Cm(G)\otimes \Z[U]$, where all~$U_i$ act on $\Z[U]$ by
  multiplication by~$U$.  We think of $C'(G)$ as a bigraded complex
  (rather than as a filtered one), writing
  $C'(G)=\bigoplus_{s\in\Z}C'(G,s)$, where here $s$ refers to the
  Alexander grading.  The homology of $C'(G)$ is easily seen to agree
  with $\HFm(G)\otimes V^{n-1}$ (following
  Lemma~\ref{lemma:TauTilde}), where the tensor product is taken in
  the bigraded sense and, as before, $V$ is a rank two module
  generated by two elements, one with Maslov and Alexander bigrading
  $(0,0)$, and another with Maslov and Alexander bigrading $(-1,-1)$.
  It follows readily that
  $$\tau'(K)=\max\,\{\,s\in\Z
  \mid\exists \xi\in H_*(C'(G,s))~\text{such that}~\forall d\geq
  0,~U^d\xi\neq0\,\}.$$

  Given any $s\in\Z$, $C'(G,s)$
  is a chain complex over $\Z$ which still retains its Maslov grading.
  It is generated by
  elements $U^{m} \cm \x$ with $m\geq 0$ and $A(\x)-m=s$, its
  differentials count those empty rectangles with $\sum_i X_i(r)=0$,
  and each rectangle is counted with multiplicity $U^{\sum_i O_i(r)}$.
  Let $C''(G,s)$ be the chain complex generated over $\Z$ by those
  $\x\in\S$ with $A(\x)\geq s$ and differential
  \begin{equation}
    \label{eq:PartialPP}
    \partial''(\x) = \sum_{\y\in\S(G)}\,\sum_{\substack{r\in\EmptyRect(\x,\y)\\r \cap \Xs = \emptyset}}
\sign(r)\cm \y.
\end{equation}
  There is a canonical inclusion of complexes $\iota\colon C''(G,s) \subset C''(G,s-1)$.
  Let $C''(G)$ be the union of all $C''(G,s)$.
  For all $s\in\Z$, there is an isomorphism of chain complexes of $\Z$-modules
  $$\phi_s \colon C'(G,s) \longrightarrow C''(G,s)$$
  defined by 
  $$\phi_s(U^m \x) = \x,$$
  which  fits
  into the diagram
  $$\begin{CD}
    C'(G,s) @>>> C''(G,s) \\
    @V{U}VV @V{\iota}VV \\
    C'(G,s-1) @>>> C''(G,s-1).
  \end{CD}$$
  Therefore the inclusion of $H_*(C'(G,s))$ into the direct
  limit of $H_*(C'(G,*))$ (with connecting map~$U$)
  corresponds to the inclusion of $H_*(C''(G,s))$ into the homology of
  $C''(G)$.

  Let $G'$ be the grid diagram obtained by reversing the roles of the
  $O_i$ and the $X_i$ in~$G$; $G'$~is a diagram for~$-K$, the knot~$K$
  with the orientation reversed.  Let $A(\x;G')$ be the Alexander
  grading with respect to $G'$.  By
  Equation~(\ref{eq:AlexanderFormulaOne}) we see that $A(\x;G') =
  -A(\x;G) -n+1$, so $C''(G,s) \cong \Filt(G',-s-n+1)$.  Therefore the
  filtered
  chain homotopy type of $C''(G)$ is identified with the filtered
  chain homotopy type of~${\widetilde C}(G')$, with modified
  filtration degree.  Thus $\tau'(K) =
  -\widetilde \tau(-K)-n+1$ and, by Lemma~\ref{lemma:TauTilde},
  $\tau'(K) = -\tau(-K)$.  We also have $\tau(K) = \tau(-K) =
  -\tau(\Mirror(K))$~\cite{4BallGenus}. (This last step can also be proved
  using grid diagrams alone,
  cf.~\cite[Proposition~\ref{OnComb:prop:Mirror}]{MOST}.)
  Putting these together, we complete the proof.
\end{proof}


\bibliographystyle{halpha}
\bibliography{biblio}

\end{document}